\documentclass[opre,nonblindrev]{informs3}

\DoubleSpacedXI 

\usepackage{endnotes}
\let\footnote=\endnote
\let\enotesize=\normalsize
\def\notesname{Endnotes}%
\def\makeenmark{\hbox to1.275em{\theenmark.\enskip\hss}}
\def\enoteformat{\rightskip0pt\leftskip0pt\parindent=1.275em
  \leavevmode\llap{\makeenmark}}

\usepackage{natbib}
 \bibpunct[, ]{(}{)}{,}{a}{}{,}%
 \def\bibfont{\small}%
\usepackage{amsmath,amsfonts,amssymb,amscd}
\usepackage{booktabs}
\usepackage{enumerate}
\usepackage{graphicx, subcaption, caption}
\usepackage[ruled,vlined]{algorithm2e}
\usepackage{hyperref}
\usepackage{nicefrac}
\usepackage{tikz}
\usetikzlibrary{decorations.pathreplacing,calligraphy, arrows.meta, spy}
\usepackage{pgfplots}
\pgfplotsset{width=8cm, height=6cm,compat=1.9}
\SetKwInput{KwInput}{Input}
\SetKwInput{KwOutput}{Output}
\SetKwInput{KwInitialize}{Initialize}

\newcommand{\ud}{\,\mathrm{d}}
\newcommand{\q}{\frac}

\newcommand{\EE}{{\mathbb{E}}}
\newcommand{\e}{{\mathsf{e}}}
\newcommand{\F}{{\mathsf{F}}}

\newcommand{\N}{{\mathbb{N}}}
\renewcommand{\P}{{\mathcal{P}}}
\newcommand{\PP}{{\mathbb{P}}}
\newcommand{\R}{{\mathcal{R}}}
\newcommand{\RR}{{\mathbb{R}}}

\newcommand{\Var}{{\mathsf{Var}}}
\newcommand{\x}{{\boldsymbol{x}}}
\newcommand{\X}{{\mathcal{X}}}
\newcommand{\opt}{{\mathsf{OPT}}}
\newcommand{\alg}{{\mathsf{ALG}}}
\newcommand{\flu}{{\mathsf{FLU}}}
\newcommand{\birt}{{\mathsf{BIRT}}}

\newcommand{\norm}[1]{\left\lVert#1\right\rVert}

\newenvironment{shawn}{\color{red}}{}

\TheoremsNumberedThrough     
\ECRepeatTheorems

\EquationsNumberedThrough    

\begin{document}


\RUNAUTHOR{Balseiro and Xia}

\RUNTITLE{Uniformly Bounded Regret in Dynamic Fair Allocation}
\TITLE{Uniformly Bounded Regret in \\Dynamic Fair Allocation}

\ARTICLEAUTHORS{%
\AUTHOR{Santiago R. Balseiro, Shangzhou Xia}
\AFF{}
\AFF{Graduate School of Business, Columbia University, New York, NY 10027, \EMAIL{srb2155@columbia.edu}, \EMAIL{sx2182@columbia.edu}}
} 

\ABSTRACT{%
We study a dynamic allocation problem in which $T$ sequentially arriving divisible resources are to be allocated to a number of agents with linear utilities. The marginal utilities of each resource to the agents are drawn stochastically from a known joint distribution, independently and identically across time, and the central planner makes immediate and irrevocable allocation decisions. Most works on dynamic resource allocation aim to maximize the utilitarian welfare, i.e., the efficiency of the allocation, which may result in unfair concentration of resources on certain high-utility agents while leaving others' demands under-fulfilled. In this paper, aiming at balancing efficiency and fairness, we instead consider a broad collection of welfare metrics, the H\"older means, which includes the Nash social welfare and the egalitarian welfare. To this end, we first study a fluid-based policy derived from a deterministic surrogate to the underlying problem and show that for all smooth H\"older mean welfare metrics it attains an $O(1)$ regret over the time horizon length $T$ against the hindsight optimum, i.e., the optimal welfare if all utilities were known in advance of deciding on allocations. However, when evaluated under the non-smooth egalitarian welfare, the fluid-based policy attains a regret of order $\Theta(\sqrt{T})$. We then propose a new policy built thereupon, called Backward Infrequent Re-solving with Thresholding ($\mathsf{BIRT}$), which consists of re-solving the deterministic surrogate problem at most $O(\log\log T)$ times. We prove the $\mathsf{BIRT}$ policy attains an $O(1)$ regret against the hindsight optimal egalitarian welfare, independently of the time horizon length $T$. We conclude by presenting numerical experiments to corroborate our theoretical claims and to illustrate the significant performance improvement against several benchmark policies.
}%



\maketitle

%

\section{Introduction}
\label{sec:intro}
Fair resource allocation has been widely studied in economics and operations research, with applications in school seat allocation~\citep{AbdulkadirougluSo03, AbdulkadirogluChPaRo17}, work visa lotteries~\citep{PathakReSo20}, affordable housing allocation~\citep{Kaplan84, ArnostiSh20}, machine load balancing~\citep{Azar95, Azar98}, etc. Involving multiple parties stipulates fairness as a crucial concern to the central planner. To this end, allocations are often evaluated on a balance between efficiency and fairness across the agents~\citep{BertsimasFaTr11, BertsimasFaTr12, AhmedDiFu17, DonahueKl20, ZengPs20}. Common fairness criteria studied in the literature of fair allocation include Pareto efficiency~\citep{Varian73}, envy-freeness~\citep{AbdulkadirogluChPaRo17} along with its variants, and proportional fairness~\citep{BrandtCoEn16, BanerjeeGkHo22}. 

Cardinal welfare metrics incorporate efficiency and fairness concerns of allocations by mapping agents' individual utilities to a real-valued numeric representation of collective welfare. These metrics are wieldy in operations research and computer science because they allow to cast allocation problems as optimization problems, and different allocations may be compared against one another, yielding a quantifiable welfare difference. One particular cardinal fairness metric is the egalitarian welfare, based on the theory of distributive justice by~\citet{Rawls01, Rawls20}, whose difference principle maximizes the welfare of the worst-off group of agents. The egalitarian welfare objective defines a max-min welfare optimization problem, whose indivisible variant is dubbed the Santa Claus problem~\citep{BansalSv06, Ben98, ChakrabartyChKh09}. Another well-studied cardinal example is the Nash social welfare, mathematically the geometric mean of the agents' utilities~\citep{Nash50, LuceRa89, HahnArIn82, FreemanAhVa20, CaragiannisKuMoPr19, ConitzerFrSh17}. Optimal allocations under the Nash social welfare are invariant to scaling of any single agent, for it is homogeneous on each agent's utility. In this work we consider a broad parameterized collection of \emph{H\"older mean} welfare metrics~\citep{WuLiGa21}, which subsumes the two examples above.

Allocation problems in real life, however, often concern dynamic and irrevocable decisions where either agents or resources arrive in a sequential manner~\citep{NandaXuSaDi20, DickersonSaSrXu19, CheungLyTeWa20}. Some literature has referred to this variant as the online rationing problem, with applications in allocating computational resources to cloud users~\citep{DinhLiQuSh20, ZhangXiZhLi18, CayciGuEr20} and most recently in medical equipment and vaccine rationing during the COVID-19 pandemic~\citep{DevanurJaSiWi19, Grigoryan21, AlkaabnehDiGa21}.

This work focuses on the specific variant of allocating resources arriving sequentially to heterogeneous agents, which is particularly relevant when resources are perishable or urgently needed and require immediate allocation decisions. Applications include allocating donated organs to hospitals~\citep{AgarwalAsReSo19}, donated food to local charities~\citep{AleksandrovAzGaWa15, LienIrSm14}, and online advertisement slots to advertisers~\citep{ChoiMeBaLe20, Watts21, GollapudiPa14, LiRoZhZh21, BateniChCiMi22, ArnostiBeMi16}, etc. Our objective is to design computationally efficient dynamic policies that maximize the overall welfare.

\paragraph{Contributions} In this work, we consider the dynamic fair allocation problem, where $T$ resources arriving sequentially are to be immediately and irrevocably allocated to $n$ agents. The joint marginal utilities of each resource to the agents are drawn independently from a joint distribution that is known to the central planner, which are revealed before allocation decisions are made. An optimal online policy could in principle be computed using dynamic programming but, because of the so-called \textit{curse of dimensionality}, solving the problem to optimality is impractical when the time horizon is large. We therefore seek to design policies that are computationally efficient and have provable performance guarantees.

We first consider a static policy that is obtained by solving a fluid relaxation in which all stochastic quantities are replaced by their mean values. Simple as it is, we show under all smooth H\"older-mean welfare metrics it attains a resounding $O(1)$ regret over any $T$ against the hindsight optimal solution, i.e., the optimal allocation if all marginal utilities were known in advance. The sole non-smooth exception is the egalitarian welfare objective, under which fluid policies attain a regret on the order of $\Theta(\sqrt{T})$. Refer to Section~\ref{sec:fluid} for rigorous definitions of the fluid relaxation problem, its optimum $\flu$, and fluid policies $\mathsf{F}$. 

Our next main contribution is proposing the Backward Infrequent Re-solving with Thresholding ($\birt$) policy for dynamic fair allocation under the egalitarian welfare. Our policy is inspired by~\citet{BumpensantiWa20}, who designed similar policies for utilitarian objectives based on fluid policies, and consists of updating and re-solving the fluid problem at most $O(\log\log T)$ times. Utilitarian objectives are simpler to handle because they are separable over time, i.e., the total utility is simply the sum of the utilities of each time period. By contrast, the egalitarian objective discussed in our paper is non-linear and not time separable, and so requires a different analysis. We provide a novel analysis to show that the $\birt$ policy achieves $O(1)$ regret against the hindsight optimum, namely uniformly bounded over the time horizon length $T$. The uniformly bounded regret guarantee implies that the $\birt$ policy performs almost comparably to a clairvoyant who foresees all arrivals in advance of acting.

A salient feature of the $\birt$ policy is that its $O(1)$ regret guarantee is insensitive to whether or not the original fluid problem is degenerate (more technically, that any fluid policy is nondegenerate) or nearly degenerate. This is in drastic contrast to previous works on online stochastic optimization whose similar $O(1)$-loss results rely heavily on non-degeneracy conditions on the underlying fluid problem (see Section~\ref{sec:experiment} for a numerical illustration). 

\begin{figure}[htbp]
\centering
\small
\begin{tikzpicture}[ultra thick, scale=0.48]
\draw [->] (0,-4) -- (0,4);
\filldraw [white] (0,4.65) circle (0pt) 
node[anchor=south, black]{Expected};
\filldraw [white] (0,4) circle (0pt) 
node[anchor=south, black]{welfare};
\filldraw [white] (-5,4.75) circle (0pt) 
node[anchor=south, black]{\underline{$q=-\infty$}};
\filldraw [white] (-5,5.2) circle (0pt) 
node[anchor=north, black]{egalitarian};
\filldraw [white] (5,4.75) circle (0pt) 
node[anchor=south, black]{\underline{$q\in(-\infty,1]$}};
\filldraw [white] (5,5.2) circle (0pt) 
node[anchor=north, black]{H\"older-mean};

\draw [orange] (-0.2,3) -- (0.2,3); 
\filldraw [white] (0.2,3) circle (0pt)
node[anchor=west, orange]{$\mathsf{FLU}$};

\draw [orange] (-0.2,0) -- (0.2,0); 
\draw [line width=0.1mm] (-1.5,0) -- (-1.1,0); 
\filldraw [white] (0.2,0) circle (0pt)
node[anchor=west, orange]{$\mathsf{OPT}$};

\draw [blue] (-0.2,-1) -- (0.2,-1); 
\draw [line width=0.1mm] (-1.5,-1) -- (-1.1,-1); 
\filldraw [white] (0.2,-1) circle (0pt)
node[anchor=west, blue]{$\mathsf{BIRT}$};

\draw [blue] (-0.2,-3) -- (0.2,-3); 
\filldraw [white] (0.2,-3) circle (0pt)
node[anchor=west, blue]{$\mathsf{F}$};

\draw [decorate,
    decoration = {
    calligraphic brace,
    amplitude=6pt}] (-0.4,0.05) --  (-0.4, 2.95);
\filldraw [white] (-1.5,1.5) circle (0pt)
node[anchor=east, black]{[Lemmata~\ref{lemma:opt<=flu}, \ref{lemma:flu-opt}] $\Theta(\sqrt{T})$};

\draw [decorate,
    decoration = {
    calligraphic brace,
    amplitude=6pt}] (-0.4,-2.95) --  (-0.4,-0.05);
\filldraw [white] (-1.5,-1.5) circle (0pt)
node[anchor=east, black]{[Lemmata~\ref{prop:regret-fluid-sqrt},~\ref{lemma:opt-flu}] $\Theta(\sqrt{T})$};

\draw [<->, line width=0.2mm] (-1.3,0) -- (-1.3,-1); 
\filldraw [white] (-2,-0.5) circle (0pt)
node[anchor=east, red]{[Theorem~\ref{thm:main}] $\; O(1)$};

\draw [decorate, decoration = {calligraphic brace, amplitude=6pt}] (2.3,2.95) --  (2.3,-2.95);
\filldraw [white] (2.7,0) circle (0pt)
node[anchor=west, red]{$O(1)$ [Theorem~\ref{thm:main-fluid}]};

\end{tikzpicture}
\caption{Summary of benchmark values and policy performances. Results relevant to the egalitarian welfare are shown on the left, and that relevant to all other H\"older-mean welfare metrics is shown on the right.}
\label{fig:axis}
\end{figure}
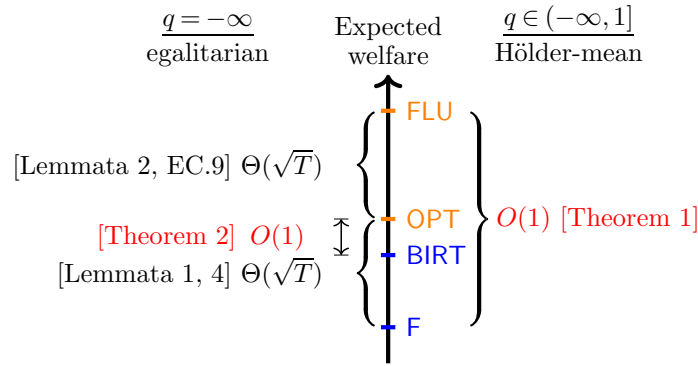

We summarize findings on the benchmarks ($\flu$ and $\opt$) and policy performances ($\birt$ and $\F$) on the vertical axis in Figure~\ref{fig:axis} with respective suprema of the differences (i.e., worst-case arrival distribution) between them. All in all, we provide computationally efficient online policies that attain uniformly bounded regret in for all H\"older-mean welfare metrics.

\subsection{Related work} 
Our work contributes to the literature on fair resource allocation and, more generally, to the literature on dynamic resource allocation problems.

\paragraph{Fair allocation} There is a stream of literature studying the online version of the problem in which the arrival process in unknown to the central planner. \citet{DevanurJaSiWi19} consider the egalitarian objective when the agents' utilities are drawn independently from a distribution that is unknown to the central planner (the ``unknown IID'' model) and design an algorithm that attains $O(\sqrt{T})$ regret. \citet{AgrawalDe14} consider general online allocation problems with concave non-separable objectives, which subsume our problem with H\"older-mean welfare metrics, and propose an algorithm based on multiplicative weight updates method~\citep{AroraHaKa12} that also attains $O(\sqrt{T})$ regret. \citet{BalseiroLuMi21} study similar general online allocation with hard resource constraints and propose fast algorithms, which do not require solving auxiliary optimization problems, based on mirror descent. More recently, \citet{KawaseSu21} consider the max-min fair allocation of indivisible resources arriving online to heterogeneous agents. They propose a deterministic policy that achieves an asymptotic approximation ratio of $n$ in the adversarial arrival model, which they show is optimal. In the case of unknown IID arrivals, they give a policy that achieves $O(\sqrt{T})$ regret against the hindsight optimum benchmark, based on the multiplicative weights updates~\citep{AroraHaKa12}.  

In contrast to the previous line of work, we consider the case where the distribution of utilities is known to the central planner. The $O(1)$ regret bound of the $\birt$ policy in the known IID model significantly improves upon the unknown IID model, in which the best known attainable regret is $O(\sqrt{T})$~\citep{AgrawalDe14, DevanurJaSiWi19, KawaseSu21, BalseiroLuMi21}. The drastic performance improvement shows the great use of distributional knowledge. This is similar in spirit to~\citet{BanerjeeGkGoJi22}, who show in maximizing Nash social welfare, agents' utility predictions significantly improve the approximation ratio of online algorithms in adversarial settings to $O(\log (n \wedge T))$ from a trivial $\Theta (n)$. Other works on maximizing Nash social welfare include~\citet{FreemanSeVi17, SinclairJaBaYu20}.

Under the egalitarian welfare, the particular case of identical agents has been studied for decades as the ``online machine covering'' problem in the context of scheduling~\citep{DeuermeyerFrLa82, AzarEp98, GalvezSoVe20, Woeginger97}. Moreover, there is a line of work studying an offline version of the max-min fair allocation problem (the case where all utilities are known in advance), which is referred to as the Santa Claus problem~\citep{LiHeJiWu19, AsadpourSa10, KawaseSu20, BezakovaDa05, Golovin05, Feige08, ChakrabartyChKh09, HaeuplerSaSr11}.

In this paper, we study an \emph{online supply} version of the dynamic fair allocation problem in which resources arrive online and the population of agents are fixed. The \emph{online demand} version of the problem in which resources are fixed and agents arrive online has also been studied in the literature~\citep{NandaXuSaDi20, ManshadiNiRo21, MaXuXu20, LienIrSm14}.

Finally, there is a long line of literature on online allocation focusing on other fairness notions such as envy~\citep{ZengPs20, BenadeKaPrPs18}, envy-freeness (up to 1 item)~\citep{HePrPsZe19, AleksandrovAzGaWa15}, Pareto efficiency~\citep{ZengPs20} and more general concave returns~\citep{DevanurJa12}. We recommend to the reader the recent survey by~\citet{AleksandrovWa20} for an overview of the literature of variants of online fair allocation problems.

\paragraph{Dynamic resource allocation}
Our work is also closely related methodologically to the stream of literature on dynamic resource allocation, especially revenue management, that aims at maximizing revenue or efficiency of an allocation given limited resources. Our policy is based on re-solving the fluid problem, a specific form of certainty equivalent of dynamic stochastic optimization problems when the problem instances are independent and identically distributed according to a distribution known to the central planner. The main difference with this line of work is that the H\"older-mean welfare objective we consider is not time separable, which calls for a different theoretical analysis.

The seminal works by~\citet{GallegoVa97} and~\citet{TalluriVa98} show fluid policies based on solving the fluid problem only once at the beginning suffice to achieve an $O(\sqrt{T})$ regret against the hindsight optimum. \citet{ReimanWa08} propose a policy based on re-solving the fluid policy exactly once only incurring an $o(\sqrt{T})$ asymptotic revenue loss. Re-solving multiple times or even every period can give an $O(1)$ loss under the assumption that the fluid problem admits a nondegenerate optimal solution~\citep{JasinKu12, WuSrLiJi15}. \citet{BumpensantiWa20} propose the Infrequent Re-solving with Thresholding policy, which incurs an $O(1)$ loss, and dispenses with the nondegeneracy condition. Our policy is based thereupon. Most recently, \citet{VeraBa21} and \citet{VeraBaGu21} propose a general Bayes selector policy that attains low regret for a number of dynamic allocation problems. A specific variant called the multi-secretary problem also witnessed policies incurring $O(1)$~\citep{ArlottoGu19} or $O(\log T)$~\citep{Bray19} regret. We recommend to the reader a recent survey by~\citet{BalseiroBePi21} for an overview of dynamic allocation problems. 

\subsection{Notation} 
We denote the extended real line by $\overline\RR := \RR \cup \{-\infty, +\infty\}$. For $m\in\N$, we denote $[m]:=\{1,2,\ldots,m\}$, the probability $m$-simplex by $\Delta_m:=\{x\in\RR_+^m:\norm{x}_1 = 1\}$, a column $m$-vector of ones by $\textbf{1}_m$ and a column $m$-vector of zeros by $\textbf{0}_m$
. For $w\in\RR$, we denote by $\lfloor w \rfloor:=\max\{m\in\mathbb{Z}:m\le w\}$ the largest integer smaller than or equal to $w$, and by $\lceil w \rceil:=\min\{m\in\mathbb{Z}:m\ge w\}$ the smallest integer larger than or equal to $w$. 
For $u,v \in \RR^n$, we denote the element-wise product by $u*v := (u^iv^i)_i$. We use boldface alphabets (e.g., $\boldsymbol\xi, \boldsymbol{x}$) to denote matrices. For any two positive-valued functions $f,g:\mathbb{N}\to\RR_+$, we denote $f(T) = o(g(T))$ if $\lim_{T\to\infty} f(T)/g(T) = 0$; we denote $f(T) = O(g(T))$ if $\limsup_{T\to\infty} f(T)/g(T) < \infty$; we denote $f(T) = \Omega(g(T))$ if $\liminf_{T\to\infty} f(T)/g(T) > 0$; we denote $f(T) = \Theta(g(T))$ if $f(T) = \Omega(g(T))$ and $f(T) = O(g(T))$. With a slight abuse of notation, we sometimes denote optimization problems and their optima interchangeably.

\section{Problem formulation}
\label{sec:model}
In this work, we consider the problem of dynamically allocating resources to $n$ heterogeneous agents with linear utilities that are stochastically distributed. More formally, suppose each agent $i\in[n]$ starts with zero initial utility $B_0^i=0$ and passively accumulates linear utility while receiving allocated resources over time. The central planner receives a sequence of $T$ divisible resources, each along with a revealed vector $b_t\in\RR_+^n$ indicating its joint marginal utilities to the $n$ agents, and immediately decides on its proportional allocation $x_t\in\Delta_n$ to the agents. In the sequel, we refer to the $t$-th resource and the $t$-th period interchangeably, and to resources and arrivals interchangeably. By the end of the time horizon, agent $i$ would have a cumulative utility of $B_T^i = \sum_{t=1}^T b_t^i x_t^i$.

We evaluate allocation decisions on their balance between efficiency and fairness through a cardinal welfare metric $w: \RR_+^n \to \RR_+$ of the vector of agents' cumulative utilities $B_T\in\RR_+^n$. We denote the welfare generated by a sequence of allocation decisions $\x\equiv (x_t:t\in[T])$ according to policy $\pi$ by 
\begin{equation}
\label{eqn:alg}
    \alg(\pi) := w \left( \sum_{t=1}^T b_t^i x_t^i \right).
\end{equation}
An \emph{online policy} is non-anticipative, or in other words, decides on allocations only based on the history observed so far. More formally, an online policy allocates $x_t \in \sigma(b_1,\ldots,b_t)$ for all $t\in[T]$, where $\sigma(b_1,\ldots,b_t)$ is the $\sigma$-algebra generated by the sequence of utility vectors. Notice we assume that online policies are deterministic for the purpose of evaluating the  welfare they generated---this is without loss of optimality because deterministic online policies weakly dominate their randomized counterparts (Lemma~\ref{lemma:futility}).

A conventional benchmark used to evaluate online policies is the \textit{hindsight optimum}, defined as the maximum welfare~\eqref{eqn:alg} if allocations $\x$ were allowed to be anticipative, i.e., all vectors of joint marginal utilities $(b_t:t\in[T])$ were known in advance. We denote by $\opt$ the hindsight optimum, as well as any such optimal policy with a slight abuse of notation. Formally, the hindsight optimum is the optimal value of the maximization problem 
\begin{equation}
\label{eqn:opt}
    \opt := \max \left\{ w \left( \sum_{t=1}^T b_t^i x_t^i \right) : \x \in \Delta_n^T \right\}.
\end{equation}

\subsection{Welfare metric}
In this work, we evaluate allocations using cardinal welfare metrics $w$, i.e., real-valued functions of the agents' cumulative utilities $B_T \in \RR_+^n$, which allow to balance between efficiency and fairness concerns. Hence, allocation decisions can be directly compared against one another, and their difference indicates a quantifiable welfare loss. 

We next enumerate three cardinal welfare examples and briefly discuss the intuition behind them. The \emph{utilitarian welfare} $\sum_{i\in[n]} B_T^i$ neglects any fairness concerns; optimizing it leads to efficient allocations that could and most often lead to highly unfair concentration of resources at certain agents or groups. The \emph{egalitarian welfare} $\min_{i\in[n]} B_T^i$ is the other extreme---defined as the utility of the worst-off agent, this metric provides all agents with a uniform minimal utility guarantee. The \emph{Nash social welfare} $(\prod_{i=1}^n B_T^i)^{1/n}$ is a more balanced metric in terms of the trade-off between efficiency and fairness; in addition, its geometric mean form satisfies positive homogeneity in \emph{each agent's utility}, which insulates the central planner against agents' possible strategic reporting of marginal utilities.

Formally, we consider a collection of welfare metrics that subsumes and generalizes the examples above, the \emph{H\"older means} $w_q:\RR_+^n \to \RR_+$, a.k.a.\ \emph{generalized power means}, parameterized by $q\in\overline{\RR}$, defined as follows and illustrated in the fairness spectrum shown in Figure~\ref{fig:spectrum}. In this problem, we stipulate $q\in[-\infty,1]$ because $w_q$ for $q\in(1,+\infty]$ is not concave and hence deemed ``unfair'' as per Lemma~\ref{lemma:axioms}.
\begin{itemize}
\item For $q=-\infty$, $w_q(B) = \min_i B^i$;
\item for $q\in(-\infty,0)$, $w_q(B) = (\q1n\sum_{i\in[n]} (B^i)^q)^{1/q}$ if $B>0$ and $w_q(B)=0$ otherwise; 
\item for $q=0$, $w_q(B) = (\prod_{i=1}^n B_i)^{1/n}$; 
\item for $q\in(0, +\infty)$, $w_q(B) = (\q1n\sum_{i\in[n]} (B^i)^q)^{1/q}$; 
\item for $q=+\infty$, $w_q(B) = \max_i B^i$.
\end{itemize}
{
\begin{figure}[th]
\centering
\begin{tikzpicture}[scale=0.9]
\draw [<-,ultra thick, blue] (-5,0) -- (2,0);
\draw [->,ultra thick, blue,dotted] (2,0) -- (5,0);
\filldraw [black] (-6, 0.5) circle (0pt) node[anchor=east]{parameter $q$};
\filldraw [blue] (-5.1, 0) circle (0pt) node[anchor=east]{fairer};
\filldraw [black] (-6.5, -0.4) circle (0pt) node[anchor=east]{H\"older-mean};
\filldraw [black] (-5.9, -0.6) circle (0pt) node[anchor=east]{$w_q$};
\filldraw [black] (-6.5, -0.9) circle (0pt) node[anchor=east]{welfare metric};
\filldraw [blue] (5.1, 0) circle (0pt) node[anchor=west]{more unfair};

\draw [blue,thick] (-5, 0.2) -- (-5, -0.2); 
\filldraw  (-5, 0.2) circle (0pt)
node[anchor=south]{$-\infty$};
\filldraw [black] (-4.6, -0.2) circle (0pt)
node[anchor=north]{egalitarian};
\filldraw [black] (-4.6, -0.7) circle (0pt)
node[anchor=north]{(minimum)};

\draw [blue,thick] (-2, 0.2) -- (-2, -0.2); 
\filldraw  (-2, 0.2) circle (0pt)
node[anchor=south]{$-1$};
\filldraw [black] (-2, -0.2) circle (0pt)
node[anchor=north]{harmonic};
\filldraw [black] (-2, -0.7) circle (0pt)
node[anchor=north]{mean};

\draw [blue,thick] (-0, 0.2) -- (-0, -0.2); 
\filldraw  (0, 0.2) circle (0pt)
node[anchor=south]{$0$};
\filldraw [black] (0, -0.3) circle (0pt)
node[anchor=north]{Nash};

\draw [blue,thick] (2, 0.2) -- (2, -0.2); 
\filldraw  (2, 0.2) circle (0pt)
node[anchor=south]{$1$};
\filldraw [black] (2, -0.3) circle (0pt)
node[anchor=north]{utilitarian};

\draw [blue,thick] (5, 0.2) -- (5, -0.2); 
\filldraw  (5, 0.2) circle (0pt)
node[anchor=south]{$+\infty$};
\filldraw [black] (5, -0.3) circle (0pt)
node[anchor=north]{maximum};
\end{tikzpicture}
\caption{Fairness spectrum of the H\"older-mean welfare metrics $w_q$ parameterized by $q\in\overline\RR$. The fairness of the welfare metric decreases as $q$ increases. Dashed line indicates the ``unfair'' regime ($q>1$).}
\label{fig:spectrum}
\end{figure}%
}
We consider the broad collection of H\"older mean welfare metrics $w_q$ because they satisfy a number of axioms that are reasonable in common scenarios of fair allocation; moreover, classical results have shown that they are the only welfare metrics satisfying these axioms, up to a multiplicative constant, which we have disposed of through normalization of $w_q$ in the definition above. Hence, $w_q(u\textbf{1})=u$ for any scalar $u\ge0$, i,e., the welfare equals each agent's utility if all utilities are identical.
\begin{lemma}[\cite{Moulin04, HardyLiPo52}]
\label{lemma:axioms}
For any welfare metric $w:\RR_+^n \to \RR$, $w \propto w_q$ for some $q\in[-\infty, 1]$ if and only if $w$ satisfies the following axioms.
\begin{enumerate}
\item Monotonicity: $w(B)\ge w(\tilde B)$ if $B\ge \tilde B$.
\item Symmetry: $w(B) = w(\sigma(B))$ for any permutation $\sigma\in\mathbb{S}_n$.
\item Continuity: $w$ is continuous.
\item Independence of unconcerned agents: if $w(B^1, B^2, \ldots, B^n) \ge w(\tilde B^1, B^2, \ldots, B^n)$, then $w(B^1, \tilde B^2, \ldots, \tilde B^n) \ge w(\tilde B^1, \tilde B^2, \ldots, \tilde B^n)$.
\item Homogeneity: $w(\lambda B) = \lambda w(B)$ for $B\ge0$ and $\lambda>0$.
\item Pigou–Dalton principle: $w(B^1, B^2, B^3,\ldots) \ge w(B^1+\nicefrac\varepsilon2, B^2-\nicefrac\varepsilon2, B^3, \ldots)$ for $\varepsilon\in(0,B^2-B^1)$ if $B^1<B^2$.
\end{enumerate}
\end{lemma}

Note the axioms serve as the foundation to fair allocation and are reasonably motivated. More specifically, monotonicity establishes that larger agents' utilities are preferable; symmetry guarantees equal treatment by making the welfare metric independent of the identities of the agents; independence of unconcerned agents determines that an agent's preferences do not change if we modify the utility of other agents; homogeneity guarantees that allocation is independents of the units used to measured the utility; the Pigou-Dalton principle stipulates that the welfare metric prefers fairness to unfairness, i.e., transferring utility from a high-utility to a low-utility individual should increase fairness. By imposing these axioms, we restrict our attention to a class of welfare metrics that are rich enough to incorporate flexibility in fairness concerns and structured enough to work with.

\subsection{Arrival model}
We consider the IID arrival model. More formally, we assume marginal utility vectors $b_t$ are distributed according to some known finite joint distribution $P$ in the class $\P$, independently across time $t$. While our model assume that the joint marginal utilities are IID across time, they may be arbitrarily correlated across agents in a given period, allowing for heterogeneity in agents to the fullest extent.

\begin{definition}
\label{def:P} 
The class of \emph{admissible} probability distributions $\P$ is the set of all finite probability distributions $P$ with support in $[0,1]^n$.
\end{definition}
To be more precise, a distribution $P$ is admissible if it has a finite support $\{ \beta_\ell \in [0,1]^{n} : \ell\in[L] \}$ and is given by $P \{b=\beta_\ell\} = p_\ell$, where $L\in\N$ is the number of distinct arrival types, and $(p_\ell)_{\ell\in[L]} \in \Delta_{L}$ is a probability vector. 
In the sequel, the underlying distribution $P$ is omitted when implied in the context. We denote the number of type-$\ell$ arrivals by $N_\ell:=\sum_{t=1}^T \textbf{1} \{b_t=\beta_\ell\}$, which then follows a multinomial distribution, i.e., $(N_\ell:\ell\in[L]) \sim \mathsf{Multinomial}(T;(p_\ell:\ell\in[L]))$ with $T$ trials and success probabilities $(p_\ell:\ell\in[L])$. 

We evaluate an online policy using its expected H\"older-mean welfare loss against the hindsight optimum, along with the relative counterpart. More rigorously, given underlying arrival distribution $P\in\P$, we define the \emph{regret} and \emph{relative regret} of a policy $\pi$ by
\begin{equation}
\label{eqn:regret}
\R_T(\pi) := \EE\left[ \opt - \alg(\pi) \right] \quad \text{ and } \quad
\rho_T(\pi) := \q{\EE\left[ \opt - \alg(\pi) \right]}{\EE[\opt]}
\end{equation}
if $\EE[\opt]>0$ and $\rho_T(\pi):=0$ otherwise.


\section{Fluid problem and static policies}
\label{sec:fluid}
\textit{Certainty equivalents} have been widely used in general dynamic stochastic optimization problems to design simple heuristics with good performance guarantees~\citep{Bertsekas12}. Certainty equivalent policies seek to approximate the offline optimization by replacing stochastic quantities with fixed and known deterministic values. A particular choice is to replace all random quantities by their means. The resulting optimization problem is usually referred to as the \textit{fluid problem}, or the deterministic program. The optimal solution to the fluid problem prescribes a simple static policy, which can be implemented in the original stochastic program.

Before moving forward, we formally define the collection of \textit{static policies}, which will play an important role in our analysis. A static policy specifies a proportional allocation vector $\xi_\ell \in\Delta_n$ over agents for each type $\ell\in[L]$ and allocates accordingly $x_t = \sum_{\ell\in[L]} \xi_\ell \textbf{1}\{b_t = \beta_\ell\}$.
With an abuse of notation, we also refer to the underlying matrix $\boldsymbol\xi \equiv (\xi_\ell:\ell\in[L])$ as a static policy. In other words, a static policy is one that allocates incoming resources only depending on their types, i.e., it makes identical allocation decisions for resources of the same type. Under static policy $\boldsymbol\xi$, the cumulative utility of agent $i$ is $B_T^i = \sum_{\ell\in[L]} N_\ell \beta_\ell^i \xi_\ell^i$. It is worth noting that any hindsight optimal policy can always be implemented as a static policy (see Lemma~\ref{lemma:opt-static} in Appendix~\ref{sec:add-results}).

We now formally define the fluid problem. Consider the IID arrival model with distribution $P\in\P$ known to the online central planner. In the fluid problem, the central planner assumes that the number of type-$\ell$ arrivals is exactly $\EE[N_\ell] = T p_\ell$. Because there is no uncertainty, it is easy to see that static policies are optimal, and the fluid problem is defined as
\begin{equation}
\label{eqn:fluid}
    \flu := {\mathrm{maximize}} \;\; \left\{ w \bigg( T\sum_{\ell\in[L]} p_\ell \beta_\ell^{i}  \xi_\ell^{i} \bigg) : \boldsymbol\xi \in \Delta_n^L \right\}.
\end{equation}
In the fluid problem, the central planner chooses an static policy that maximizes the welfare corresponding to the agents' \emph{expected} utilities, where that of agent $i$ under the static policy $\boldsymbol{\xi}$ is $\EE[B_T^i] = T\sum_{\ell\in[L]} p_\ell \beta_\ell^{i} \xi_\ell^{i}$. With an abuse of notation, we also denote by $\flu$ its optimum value, a.k.a.\ the \textit{fluid benchmark}. 

An optimal solution to the fluid problem~\eqref{eqn:fluid} is called a \textit{fluid (static) policy}, denoted by $\boldsymbol\xi^\F$ or $\F$. Fluid policies are classical heuristic policies that serve as a powerful tool for the online central planner who has no other knowledge than the underlying arrival distribution. The pseudocode for fluid policies is given in Algorithm~\ref{alg:fluid}. 

\begin{algorithm}[htbp]
\caption{Fluid Policy ($\F$)}
\label{alg:fluid}
\KwInput{time horizon length $T\in\N$, arrival distribution $P\in\P$.}
\KwInitialize{initial utilities $B_0^i \gets 0$ for $i\in[n]$;}
solve the fluid problem with initial utilities: \[\boldsymbol\xi^\F \in \arg\max \left\{ w \left( T \sum_{\ell\in[L]} p_\ell \beta_\ell^{i}  \xi_\ell^{i} \right) : \boldsymbol\xi_\ell\in\Delta_n^L \right\}\]
\For(\tcp*[h]{implement the fluid policy}){$t=1,\ldots,T$}{
    observe incoming supply with utility vector $b_t$\;
    \If{$b_t = \beta_\ell$}{
        act $x_t \gets (\xi^\F)_\ell$\;
    }
    update $B_t^i \gets B_{t-1}^i + b_t^i x_t^i$ for every agent $i \in [n]$\;
}
\KwOutput{allocations $(x_t:t=1,2,\ldots,T)$ and welfare $w(B_T)$.}
\end{algorithm}

The positive homogeneity of the welfare metric $w$ implies that fluid policies are in fact independent of the time horizon length $T$ and that the fluid benchmark $\flu$ is simply a linear function of $T$.

The classical literature in revenue management have shown that for utilitarian welfare objectives, fluid policies incur a regret at most on the order of $\sqrt{T}$~\citep{GallegoVa97, TalluriVa98}, which is tight in the worst case. Here we present similar results for the fluid benchmark and fluid policies in our dynamic fair allocation problem.

\begin{proposition}
\label{prop:regret-fluid-sqrt}
Fix $P\in\P$. Under any H\"older-mean welfare $w_q$ with $q\in[-\infty,1]$, the regret of any fluid policy $\F$ is bounded by $\R_T(\F) = O(\sqrt{T})$.
\end{proposition}
Proposition~\ref{prop:regret-fluid-sqrt} is in fact a corollary of the two following lemmata. Lemma~\ref{lemma:opt<=flu} establishes the fluid benchmark as an upper bound of the expected hindsight optimum. The result is standard and follows from Jensen's inequality and the concavity of the welfare metric (see proof in Appendix~\ref{sec:proof-opt<=flu}). Hence, the regret of any policy (e.g., a fluid policy) can be bounded by its expected welfare loss from the fluid benchmark.
\begin{lemma}
\label{lemma:opt<=flu}
For any $P\in\P$, under any H\"older-mean welfare $w_q$ with $q\in[-\infty,1]$, $\EE\left[\opt\right] \le \flu$.
\end{lemma}

The next result, Lemma~\ref{lemma:flu-flu}, establishes a stochasticity gap at most of order $\sqrt{T}$ between the welfare corresponding to the agents' expected utilities and the expected welfare under any static policy, following from a concentration of measure argument. The proof can be found in Appendix~\ref{sec:proof-flu-flu}. Proposition~\ref{prop:regret-fluid-sqrt} follows because $\flu = w(\EE[B_T])$ and $\alg(\F) = w(B_T)$ under the fluid policy $\F$.
\begin{lemma}
\label{lemma:flu-flu}
Fix $P\in\P$. Under any H\"older-mean welfare $w_q$ with $q \in [-\infty,1]$ and any static policy $\boldsymbol{\xi}$, $w_q(\EE[B_T]) - \EE[w_q(B_T)] = O(\sqrt{T})$.
\end{lemma}

Now that we have established an upper bound on the regret of fluid policies, we proceed to seek potential tightening of the bound and to present the two different performances they have in distinct regimes of the fairness spectrum (Figure~\ref{fig:spectrum}), $q\in(-\infty,1]$ and $q=-\infty$. For $q\in(-\infty,1]$, the upper bound is in fact very loose, as fluid policies are powerful enough to attain an $O(1)$ regret. In particular, the regret is exactly zero under the utilitarian welfare ($q=1$). The proof can be found in Appendix~\ref{sec:proof-main-fluid}.
\begin{theorem}
\label{thm:main-fluid}
Fix $P\in\P$. Under H\"older-mean welfare $w_q$ with $q\in(-\infty,1]$, a fluid policy $\F$ attains a regret $\R_T(\F) = O(1)$ uniformly bounded over time horizon length $T$ solely dependent on $P$.
\end{theorem}
Theorem~\ref{thm:main-fluid} shows the simple fluid static policy $\F$ can attain an impressive $O(1)$ regret. We next provide the reader with some intuition for this result. A key property of the H\"older-mean welfare $w_q$ for $q\in(-\infty,1]$ is its local strong smoothness around the expected utilities. More specifically, the time-average welfare $w_q(B_T/T)$ can be locally approximated with a quadratic form around the expected utilities $w_q(\EE[B_T]/T)$, i.e.,
\[
w_q\left(\q{B_T}T\right) - w_q\left(\q{\EE[B_T]}T\right) = \nabla w_q\left(\q{\EE[B_T]}T\right)^\top (B_T - \EE[B_T]) - O\left(\left\|\q{B_T-\EE[B_T]}T\right\|_2^2\right).
\]
In the approximate quadratic form, the first-order term has zero mean by definition, and the second-order term is $O_p(T^{-1})$ by the Central Limit Theorem, as agent $i$'s utility $B_T^i = \sum_{\ell\in[L]} N_\ell \beta_\ell^i \xi_\ell^i$ is a random walk given any static policy $\boldsymbol\xi$. Note the second-order term is negative because the welfare metric $w$ is concave. Because for some welfare metrics (such as the Nash social welfare) strong smoothness only holds locally, a technical step in the proof revolves around showing that utilities fall close their expected values with high probabilities. This implies $\EE[w(B_T)] - w(\EE[B_T]) = O(1)$, and the theorem statement follows as a result.

However, under the egalitarian welfare $w_{-\infty}$, fluid policies do not perform as well. We show that a regret of order $\sqrt{T}$ given in Proposition~\ref{prop:regret-fluid-sqrt} is essentially unimprovable for fluid policies under the egalitarian welfare $w_{-\infty}$ by providing a simple example in which the fluid policy does incur an $\Omega(\sqrt{T})$ regret. See Appendix~\ref{sec:proof-opt-flu} for the proof.
\begin{lemma}
\label{lemma:opt-flu}
There exists $P\in\P$ such that $\R_T(\F) = \Omega(\sqrt{T})$ under the egalitarian welfare $w_{-\infty}$.
\end{lemma}

\begin{figure}[hbtp]
\tikzset{new spy style/.style={spy scope={%
 magnification=5,
 size=70pt, 
 connect spies,
 every spy on node/.style={
   rectangle,
   draw,
   },
 every spy in node/.style={
   draw,
   rectangle,
   fill=gray!40,
   }
  }
 }}
\tikzset{arrow/.tip={
    Latex[length=1pt, width=1.5pt]
  }
}
\centering
\begin{subfigure}[b]{0.35\textwidth}
\centering
\begin{tikzpicture}
[new spy style, anchor=mid, x=10pt, y=20pt]
\begin{axis}[
width=\linewidth, height=\linewidth,
xmin=0, xmax=5,
ymin=0, ymax=5,
xlabel={Agent 1's utility},
ylabel={Agent 2's utility}
]
\addplot[color=blue,domain=0.1:10,samples=100]{2/x};
\addplot[color=blue,domain=0.1:10,samples=100]{4/x};
\addplot[color=blue,domain=0.1:10,samples=100]{6/x};
\addplot[color=blue,domain=0.1:10,samples=100]{9/x};
\addplot[color=blue,domain=0.1:10,samples=100]{12/x};
\addplot[color=blue,domain=0.1:10,samples=100]{16/x};
\addplot[color=blue,domain=0.1:10,samples=100]{20/x};
\draw [red] (axis cs:3, 3) circle [radius=0.18];
\draw [black, fill=black] (axis cs:3,3) circle [radius=0.02];
\draw[black,->,>=arrow](axis cs:3,3)--(axis cs:3.15,3);
\draw[black,->,>=arrow](axis cs:3,3)--(axis cs:3,2.85);
\draw[black,->,>=arrow](axis cs:3,3)--(axis cs:3,3.15);
\draw[black,->,>=arrow](axis cs:3,3)--(axis cs:2.85,3);
\coordinate (spypoint) at (axis cs:3,3);
\coordinate (magnifyglass) at (axis cs:6.6,4);
\end{axis}
\spy [blue, size=60pt] on (spypoint)
   in node[size=50pt,fill=gray!40, magnification=4] at (magnifyglass);
\end{tikzpicture}
\caption{Nash social welfare}
\end{subfigure}
\hspace{0.1\textwidth}
\begin{subfigure}[b]{0.35\textwidth}
\begin{tikzpicture}[new spy style,anchor=mid, x=20pt, y=20pt]
\begin{axis}[
width=\linewidth, height=\linewidth,
xmin=0, xmax=5,
ymin=0, ymax=5,
xlabel={Agent 1's utility},
ylabel={Agent 2's utility}
]
\draw[blue](axis cs:1,1) -- (axis cs:1,5);
\draw[blue](axis cs:1,1) -- (axis cs:5,1);
\draw[blue](axis cs:2,2) -- (axis cs:2,5);
\draw[blue](axis cs:2,2) -- (axis cs:5,2);
\draw[blue](axis cs:3,3) -- (axis cs:3,5);
\draw[blue](axis cs:3,3) -- (axis cs:5,3);
\draw[blue](axis cs:4,4) -- (axis cs:4,5);
\draw[blue](axis cs:4,4) -- (axis cs:5,4);
\draw [red] (axis cs:3, 3) circle [radius=0.18];
\draw [black, fill=black] (axis cs:3,3) circle [radius=0.02];
\draw[black,->,>=arrow](axis cs:3,3)--(axis cs:3.1,3.1);
\draw[black,->,>=arrow](axis cs:3,3)--(axis cs:3.1,2.9);
\draw[black,->,>=arrow](axis cs:3,3)--(axis cs:2.9,3.1);
\draw[black,->,>=arrow](axis cs:3,3)--(axis cs:2.9,2.9);
\coordinate (spypoint) at (axis cs:3,3);
\coordinate (magnifyglass) at (axis cs:6.6,4);
\end{axis}
\spy [blue, size=60pt] on (spypoint)
   in node[size=50pt,fill=gray!40, magnification=4] at (magnifyglass);
\end{tikzpicture}
\caption{The egalitarian welfare}
\end{subfigure}
\caption{Contour plots and infinitesimal random walks around expected utilities, illustrating corresponding welfare changes under Nash social welfare and the egalitarian welfare metrics}
\label{fig:intuition}
\end{figure}
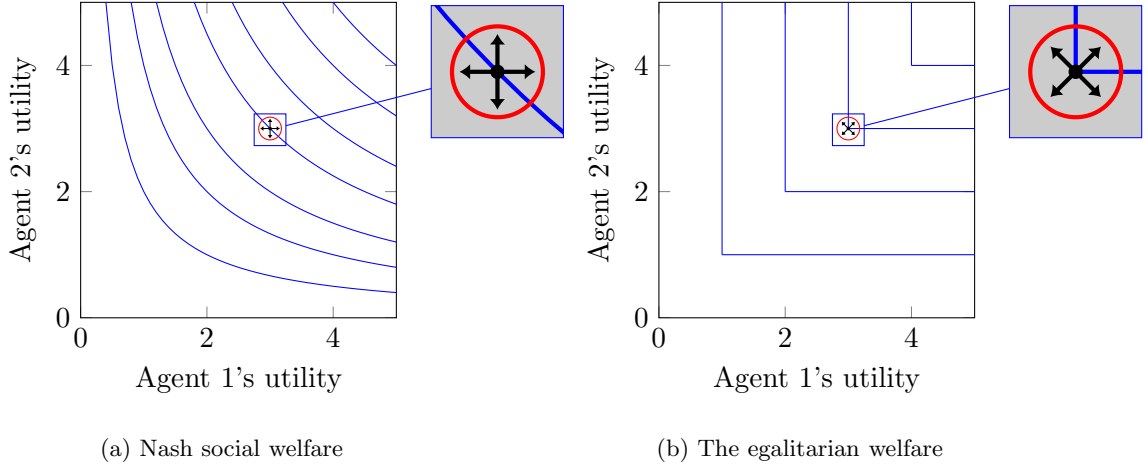

Whereas the simple fluid policies attain an $O(1)$ optimality gap for H\"older-mean welfares $w_q$ with $q\in(-\infty,1]$, they fail to attain bounded regret for the egalitarian welfare $w_{-\infty}$. In contrast to the discussion above for smooth H\"older-mean welfare metrics (i.e., $q\in(-\infty,1]$), the egalitarian welfare $\min_{i} B_T^i$ is marked by its lack of smoothness at singular points where agents have equal cumulative utilities ($B_T^1=\cdots=B_T^n$). Figure~\ref{fig:intuition} illustrates the different behavior infinitesimal perturbations around expected utilities $\EE[B_T]$ may have under smooth versus non-smooth welfare metrics. In the former, the welfare metric is locally linear, so first-order welfare changes resulting from perturbations cancel out in expectation leaving only the second-order changes; in the latter, however, the welfare metric is not differentiable, so first-order welfare changes negatively impact performance. Unfortunately, the expected utility vector always lies on a singular point, i.e., the vertex of some iso-welfare hyper-surface, because the fluid policy always generate equal expected utilities across agents (see Lemma~\ref{lemma:prophet-fair}). The fluid policy is then doomed to a first-order welfare loss in expectation and hence a regret on the order of $\Theta(\sqrt{T})$.

Therefore, the egalitarian welfare warrants a more sophisticated approach to narrowing the regret. In the next section, we propose a computationally efficient policy based on fluid policies and show it attains $O(1)$ regret under the egalitarian welfare.

\section{A policy with uniformly bounded regret under the egalitarian welfare}
\label{sec:policy}
In this section, we focus on the egalitarian welfare $w_{-\infty}$ given by $w_{-\infty}(B_T) = \min_{i\in[n]} B_T^i$. In light of fluid policies, which attain a regret of order $\Omega(\sqrt{T})$ in the worst case, we propose a collection of policies called \textbf{Backward Infrequent Re-solving with Thresholding ($\birt$)}, detailed in Algorithm~\ref{alg:birt}, inspired by previous fluid-based re-solving policies proposed in the literature, especially by~\citet{ReimanWa08} and \citet{BumpensantiWa20}. The $\birt$ policies are indexed by a hyper-parameter $\eta$ and are denoted by $\birt_\eta$, allowing for flexibility in refining the performance through control of the re-solving frequency and corresponding thresholds. The design of $\birt$ is characterized by the following three key elements: (i) the fluid problem is periodically re-solved after updating the deterministic surrogate quantities with their realizations so far, and the resulting policy is static between consecutive re-solving epochs; (ii) re-solving occurs only $O(\log\log T)$ times, based on a schedule designed such that the remaining schedule starting from each re-solving epoch is similar; and (iii) upon each solving and re-solving epoch, static policies are adjusted by a thresholding rule. 

\begin{algorithm}[htbp]
\caption{Backward Infrequent Re-solving with Thresholding ($\birt_\eta$)}
\label{alg:birt}
\KwInput{welfare metric $w:\RR_+^n\to\RR$, time horizon length $T\in\N$, hyper-parameter $\eta\in(1, \nicefrac43)$.}
\KwInitialize{set $B_0^i \gets 0$ for $i\in[n]$\;
set $K \gets \lceil \log\log T / \log \eta \rceil$\;
set $t_0^*\gets0$ and $t_k^* \gets T - \lfloor \exp(\eta^{K-k}) \rfloor$ for $k=1,\ldots,K$\;
set $\gamma_k \gets (T-t_{k+1}^*)/2n^2(T-t_k^*)$ for $k=0,\ldots,K-1$ and $\gamma_K \gets 0$\;}
\For{$k=0,1,\ldots,K$}{
    solve the updated fluid problem with updated welfare $$\boldsymbol\xi^\F \in \arg\max \left\{ w \left( B_{t_k^*} + (T-t_k^*) \sum_{\ell\in[L]} p_\ell \beta_\ell * \xi_\ell \right) : \boldsymbol\xi\in\Delta_n^L \right\};$$
    \For(\tcp*[h]{threshold the static solution}){$\ell=1,2,\ldots,L$}{
        pick some agent $j \in \arg\max_i (\xi^\F)_\ell^i$\;
        \For{$i\ne j$}{
        set $(\xi^k)_\ell^i \gets (\xi^\F)_\ell^i \, \textbf{1}\{(\xi^\F)_\ell^i\ge\gamma_k\}$\;
        }
        set $(\xi^k)_\ell^j \gets 1 - \sum_{i\ne j}(\xi^k)_\ell^i$\;
    }
    \For(\tcp*[h]{implement the thresholded static policy}){$t=t_k^*+1,\ldots,t_{k+1}^*$}{
        observe incoming supply with utility vector $b_t$\;
        \For{$\ell=1,2,\ldots,L$}{
        \If{$b_t = \beta_\ell$}{
            allocate $x_t \gets (\xi^k)_\ell$\;
        }
        }
        update welfare $B_t^i \gets B_{t-1}^i + b_t^i x_t^i$ for every agent $i \in [n]$\;
    }
}
\KwOutput{allocations $(x_t:t=1,2,\ldots,T)$ and egalitarian welfare $\min_{i\in[n]} B_T^{i}$.}
\end{algorithm}

\textbf{Updated fluid problem}. At each re-solving epoch, by updating forecast utilities with realizations observed so far, the $\birt$ policies react to the stochastic deviations from the expected path initially predicted by the fluid problem. 
More formally, with agent $i$'s cumulative utility $B_t^i = \sum_{s\in[t]} b_s^i x_s^i$ by time $t\in[T]$, where $x_s^i$ is the proportion of the $s$-th resource that is allocated to agent $i$, the updated fluid problem at time $t$ is to
\begin{equation}
\label{eqn:updated-fluid}
    \mathrm{maximize} \, \bigg\{ w \bigg( B_t^i + (T-t)\sum_{\ell\in[L]} p_\ell \beta_\ell^{i}  \xi_\ell^{i} \bigg) : \boldsymbol\xi \in \Delta_n^L\bigg\},
\end{equation}
where $\EE[B_T^i|B_t^i] = B_t^i + (T-t)\sum_{\ell\in[L]} p_\ell \beta_\ell^{i} \xi_\ell^{i}$ is the expected terminal utility of agent $i$ under static policy $\boldsymbol{\xi}$  conditional on allocations by time $t$.

\begin{figure}[htbp]
\centering
\begin{tikzpicture}
\filldraw [black] (-0.6, 0.5) circle (0pt)
node[anchor=east]{Time:};
\draw [ultra thick] (0,0) -- (10,0);
\filldraw [blue] (-0.6, -0.6) circle (0pt)
node[anchor=east]{$\birt$:};

\draw [thick] (0, 0.2) -- (0, -0.2); 
\filldraw  (0, 0.2) circle (0pt)
node[anchor=south]{$t_0^*=0$};
\filldraw [blue] (0, -0.2) circle (0pt)
node[anchor=north]{$\boldsymbol\xi^0$};

\draw [blue, thick, ->] (0.3, -0.6) -- (4.7, -0.6); 
\filldraw [blue] (2.5, -0.6) circle (0pt)
node[anchor=north]{act $\boldsymbol\xi^0$};

\draw [thick] (5, 0.2) -- (5, -0.2); 
\filldraw  (5, 0.2) circle (0pt)
node[anchor=south]{$t_1^*$};
\filldraw [blue] (5, -0.2) circle (0pt)
node[anchor=north]{$\boldsymbol\xi^1$};

\draw [blue, thick, ->] (5.3, -0.6) -- (6.7, -0.6); 
\filldraw [blue] (6, -0.6) circle (0pt)
node[anchor=north]{act $\boldsymbol\xi^1$};

\draw [thick] (7, 0.2) -- (7, -0.2); 
\filldraw  (7, 0.2) circle (0pt)
node[anchor=south]{$t_2^*$};
\filldraw [blue] (7, -0.2) circle (0pt)
node[anchor=north]{$\boldsymbol\xi^2$};

\draw [blue, thick, ->] (7.3, -0.6) -- (7.7, -0.6); 

\draw [thick] (8, 0.2) -- (8, -0.2); 
\filldraw  (8, 0.2) circle (0pt)
node[anchor=south]{$t_3^*$};
\filldraw [blue] (8, -0.2) circle (0pt)
node[anchor=north]{$\boldsymbol\xi^3$};

\filldraw  (8.8, 0.2) circle (0pt)
node[anchor=south]{$\cdots$};
\draw [blue, thick, ->] (8.2, -0.6) -- (8.5, -0.6);
\filldraw [blue] (8.5, -0.6) circle (0pt)
node[anchor=west]{$\cdots$};
\draw [blue, thick, ->] (9.1, -0.6) -- (9.4, -0.6); 

\draw [thick] (9.7, 0.2) -- (9.7, -0.2); 
\filldraw  (9.6, 0.2) circle (0pt)
node[anchor=south]{$t_{K}^*$};
\filldraw [blue] (9.7, -0.2) circle (0pt)
node[anchor=north]{$\boldsymbol\xi^K$};

\draw [blue, thick, ->] (9.8, -0.6) -- (10, -0.6); 

\draw [thick] (10, 0.2) -- (10, -0.2); 
\filldraw  (10.1, 0.25) circle (0pt)
node[anchor=south]{$T$};

\end{tikzpicture}
\caption{Re-solving schedule of the $\birt$ heuristic policies}
\label{fig:schedule}
\end{figure}
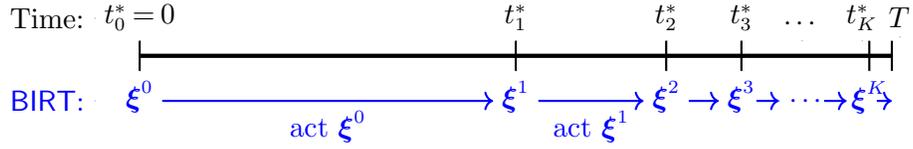

\textbf{Re-solving schedule ($\eta$)}. We now define the re-solving schedule, illustrated in Figure~\ref{fig:schedule}. Given a fixed time horizon length $T$ and a hyper-parameter $\eta\in(1, \nicefrac43)$, we set the number of re-solving epochs as $K = \lceil \log\log T/\log\eta \rceil$. The re-solving epochs, after the initial solving epoch $t_0^*=0$, are given by $t_k^* = T - \lfloor \exp(\eta^{K-k}) \rfloor$ for $k\in[K]$. Remarkably, re-solving only $O(\log\log T)$ times is  sufficient for a similar performance to that of ``frequent'' re-solving every period (see empirical evidence in Section~\ref{sec:experiment-special}). 

An important structural feature of the re-solving schedule is its self-similarity, i.e., looking into the future from any re-solving epoch gives a similar schedule. More specifically, for a re-solving epoch $\tau$ periods from the end, the next re-solving epoch is approximately $\tau^{1/\eta}$ periods from the end. We detail in Section~\ref{sec:analysis} how this greatly facilitates $\birt$'s performance analysis.

The re-solving frequency increases over time, and most re-solving epochs are concentrated near the end of the horizon, though shown spread out in Figure~\ref{fig:schedule} for illustration. 
Resolving more frequently near the end of the horizon is beneficial to the central planner because stochastic deviations from the predictions of the fluid policy are larger later on. In addition to incorporating more historical information, decisions closer to the end are more consequential due to the limited time left to make corrections. Indeed, the regret analysis in~\ref{sec:proof-main} will show that almost all regret is incurred near the end. 

\textbf{Thresholding Rule ($\gamma$)}. Given a static policy $\boldsymbol\xi\in\Delta_n^L$ and a threshold $\gamma\in(0,\nicefrac{1}{n})$, we define its corresponding thresholded static policy $\tilde{\boldsymbol\xi}$ by implementing the following procedure for each $\ell\in[L]$. After thresholding we will arrive at a static policy $\tilde{\boldsymbol\xi}$ such that $\tilde\xi_\ell^i \in \{0\} \cup [\gamma, 1]$ for all $i$ and $\ell$.  
\begin{enumerate}
    \item 
    Pick an agent that is to receive the largest type-$\ell$ allocation $j \in \arg\max_i \xi_\ell^i$.
    
    \item 
    For other agents $i\ne j$, implement allocations only if they exceed the threshold $\gamma$ and withhold the rest, i.e., set $\tilde\xi_\ell^i \gets \xi_\ell^i \, \textbf{1}\{\xi_\ell^i\ge\gamma\}$.
    
    \item
    Allocate the remaining resource to agent $j$, i.e.\ set $\tilde\xi_\ell^j \gets 1 - \sum_{i\ne j}\tilde\xi_\ell^i$.
\end{enumerate}

Thresholding is a defining characteristic of the $\birt$ policies and, as shown in the experiments in Section~\ref{sec:experiment-special}, it is crucial to obtain uniform regret guarantees. We provide some intuition about the need for thresholding in Section~\ref{sec:degeneracy}.

\subsection{Performance analysis}
\label{sec:analysis}
In this section, we first present the main result on a uniform performance guarantee of the $\birt$ policies and then provide intuition on the regret analysis. A full proof is provided in Appendix~\ref{sec:proof-main}.

\begin{theorem}
\label{thm:main}
Fix any $P\in\P$ and $\eta \in (1, \nicefrac43)$. Then, under the egalitarian welfare metric, \[\R_T(\birt_\eta) = O(1).\] The constant factor is determined solely by $P$ and $\eta$.
\end{theorem}

\noindent \textbf{Remark}. More generally, uniformly bounded regret can be attained by the $\birt_\eta$ policy, even if we do not fix the distribution $P$ but rather allow for some dependency on the time horizon length $T$. Specifically, if $\min_\ell p_\ell = \Omega (T^{-\nicefrac{1}{2} + \varepsilon})$ for some $\varepsilon\in(0,\nicefrac{1}{2}]$, then $\R_T(\birt_\eta) = O(1)$ for choices of $\eta\in(1, (1-\nicefrac{\varepsilon}{2})^{-1})$. See Theorem~\ref{thm:extension} in Appendix~\ref{sec:add-results} for a formal exposition.

The implication of Theorem~\ref{thm:main} is twofold. The regret of non-anticipative policies against the hindsight optimum can be attributed to either an intrinsic lack of clairvoyance or suboptimal policy design. The uniform bound demonstrates that our $\birt$ policies are near-optimal in both dimensions. On the one hand, their low regret gives credit to the policy design, in addition to the computational efficiency. On the other hand, with the non-anticipative $\birt$ attaining performance almost comparable to a clairvoyant, Theorem~\ref{thm:main} showcases that the central planner is able to greatly benefit from distributional knowledge.

We next give an outline of the proof, inspired by the analysis in~\cite{BumpensantiWa20}, hinging on the structure of $\birt$. In particular, the self-similar re-solving schedule warrants a serial regret analysis. More specifically, the regret of $\birt$ can be decomposed into a series of expected egalitarian welfare losses against the hindsight optimum, each attributed to one re-solving epoch in the schedule. To this end, we introduce a series of auxiliary policies $\opt^k$ for $k=0,1,\ldots,K$, where $\opt^k$ acts exactly the same as $\birt$ does until $t_k^*$ and allocates in the hindsight optimal way afterwards, i.e., it picks the optimal allocation for the remaining time periods given the decisions up to time $t_k^*$. With an abuse of notation, we also denote by $\opt^k$ the egalitarian welfare generated by the $\opt^k$ policy. More precisely, 
\begin{equation}
\label{eqn:opt-k}
    \opt^k :=
    \max\left\{ w \left( B_{t_k^*} + \sum_{t=t_k^* +1}^{T} b_t * x_t \right) : x_t \in\Delta_n, t = t_k^*+1, \ldots, T \right\},
\end{equation}
where $B_{t}^i$ is agent $i$'s welfare by $t$ under the $\birt$ policy. In particular, $\opt = \opt^0$. Notice, however, the auxiliary $\opt^k$ policy is not an online policy, as it is anticipative starting from $t_k^*+1$. 

The auxiliary policies are introduced to quantify the egalitarian welfare loss due to each re-solving epoch. For example, comparing $\opt$ and $\opt^1$ would give us that in the first epoch, as $\opt^1$ acts optimally with clairvoyance after $t_1^*$. Hence, we can decompose the regret of $\birt$ into
\begin{equation*}
\begin{split}
    \R_T(\birt) = \EE[\opt - \alg(\birt)] = \sum_{k=1}^K \EE \left[\opt^{k-1} - \opt^k \right] + \EE \left[\opt^K - \alg(\birt) \right].
\end{split}
\end{equation*}
The telescopic decomposition is useful for two reasons. First, since the problems $\opt^{k-1}$ and $\opt^{k}$ share the same decisions up to time $t_{k-1}^*$ and the re-solving schedule is self-similar, their difference in egalitarian welfare $\opt^{k-1} - \opt^k$ is isomorphic to $\opt - \opt^1$, with appropriately chosen initial welfare and horizon lengths.
Second, we can analyze each term $\EE[\opt^{k-1} - \opt^k]$ without considering inter-epoch interactions, i.e., consequences of previous allocations on subsequent periods. In other words, analyses of disjoint epochs are decoupled from one another. This follows because we can provide guarantees for each term that are \emph{uniform over the agents' initial utilities}.

The decomposition warrants an inductive argument to show the main theorem using the following result on one-epoch expected egalitarian welfare loss. For the purpose of providing a uniform framework for analyzing each re-solving epoch, we assume the agents start with initial utilities $B_0\in\RR_+^n$ that may be positive, so that the analysis can be directly deployed to subsequent epochs.

\begin{proposition}
\label{prop:one-epoch}
Fix an arrival distribution $P\in\P$, initial utilities $B_0\in\RR_+^n$ and time horizon length $T>1$. Suppose the first re-solving occurs at some $t_1\in[T]$, and the initial threshold is set as $\gamma_0 = (T-t_1)/2n^2T$, then for a constant $C>0$ determined by $P$ and independent of $T$, $ t_1$ and $B_0$,
\begin{equation}
    \EE \left[\opt - \opt^1|B_0 \right] \le
    3L t_1 \exp\left( -C \q{(T- t_1)^4}{T^3} \right).
\end{equation}
\end{proposition}
The complete proof is postponed to Appendix~\ref{sec:proof-one-epoch}. Its main idea is to focus on the coupling event $\{\opt=\opt^1\}$, where the egalitarian welfare loss is zero. On the complementary event $\{\opt\ne\opt^1\}$, we can control the difference using the almost sure uniform bound $\opt - \opt^{1} \le t_1$ (Lemma~\ref{lemma:uniform}), which follows because any loss the clairvoyant $\opt^1$ can incur against $\opt$ must only stem from the suboptimal allocations by the static policy during the first $t_1$ periods. The most complex step of the proof involves showing that $\{\opt = \opt^1\}$ is a high probability event. To do so, we argue that $\opt$ and $\opt^{1}$ generate the same egalitarian welfare if $\opt^1$ manages to replicate the exact aggregate allocations of each type to each agent made by $\opt$; we invoke sensitivity analysis results for linear programming to provide sufficient conditions in terms of the arrival process; we then prove that these sufficient conditions hold with high probability using concentration inequalities. Thresholding the fluid policy is critical to proving that the concentration inequalities hold uniformly over the initial agents' utilities as, otherwise, $\birt$ would likely fall into the trap of not being able to replicate hindsight optimal allocations due to irrevocability of allocations.

\subsection{Degeneracy and the need for thresholding}\label{sec:degeneracy}

The fair allocation problem under the egalitarian welfare can be cast as a linear program by writing the objective in the epigraph form:
\begin{equation}\label{eq:egalitarian-LP}
\begin{aligned}
\max_{u,\boldsymbol\xi} \quad\quad\quad u & \\
\text{s.t.} \quad\quad\quad  u &\le \sum_{\ell\in[L]} \EE[N_\ell] \beta_\ell^i \xi_\ell^i, \quad\quad \forall i \in [n], \\
\sum_{i\in[n]} \xi_\ell^i & \le 1, \quad\quad \forall \ell \in [L], \\
\xi_\ell^i & \ge 0, \quad\quad \forall \ell \in [L], i \in [n].
\end{aligned}
\end{equation}
We refer to the fluid problem as \emph{degenerate} if any of its optimal solutions is degenerate as defined for linear programs. More precisely, a fluid policy $\F$ is degenerate if more constraints are active than the dimension of optimization variables, i.e.,
\[\bigg|\bigg\{i\in[n] : \flu = \sum_{\ell\in[L]} \EE[N_\ell] \beta_\ell^i \xi_\ell^i \bigg\}\bigg| 
+ \bigg|\bigg\{\ell\in[L] : \sum_{i \in [n]} (\xi^\F)_\ell^i = 1 \bigg\}\bigg| 
+ \left|\left\{(\ell,i): (\xi^\F)_\ell^i = 0 \right\}\right| 
> nL + 1.\]

When the fluid problem \eqref{eq:egalitarian-LP} is non-degenerate, the same basis is optimal both in the hindsight problem \eqref{eqn:opt} when the number of arrivals of each type $N_\ell$ is close to its expected value $T p_\ell$ and in the updated fluid problem \eqref{eqn:updated-fluid} when the cumulative agents' utilities $B_t^i$ follow likely trajectories. Therefore, with high probability, the hindsight problem and the updated fluid problems have the same optimal basis, which leads to similar allocations. When the problem is non-degenerate, a vanilla resolving policy without thresholding can achieve constant regret.


Thresholding plays a key role when the fluid problem is degenerate (or close to degenerate). In this case, slight perturbations in the number of arrivals of each type or the agents' utilities can lead to changes in the optimal basis of the linear program. A vanilla resolving policy without thresholding could introduce welfare losses by choosing bases different from that of the hindsight problem. Intuitively, by adjusting small allocation components to zero, thresholding essentially forces adopting a degenerate basic solution to the fluid problem corresponding to multiple bases. 
Thresholding postpones selecting any basis by placing $\birt$ on multiple bases that are sufficiently close by so that past suboptimal allocation decisions could be made up for later once the final choice of basis is revealed. In a word, the gist of thresholding is to \emph{postpone choosing any of neighboring bases until subsequent arrivals reveal more information}.

An alternative way to understand thresholding is through our theoretical regret analysis in Section~\ref{sec:analysis}. The intuition behind thresholding is that $\birt$ aims to recover the hindsight optimum by the end of the horizon; in particular and more ambitiously, it aspires to replicate the exact amount of resources of each type $\ell$ allocated to each agent $i$ by the hindsight optimal policy. Though re-solving helps approximate these hindsight optimal allocations, exact replication is hindered by stochasticity, especially for $(\ell,i)$ where the fluid policy yields allocations $(\xi^\F)_\ell^i$ close to zero. As a result of stochasticity, the hindsight optimal policy may possibly end up allocating none of type-$\ell$ resources to agent $i$, in which case the replication goal fails completely as allocations are irrevocable. Hence, it is best to altogether adjust these allocations to zero, leaving more leeway for better informed allocations in the future. 

Moreover, thresholding becomes more aggressive with larger thresholds over time. While thresholding allows to postpone allocation decisions that are close to zero or one, it introduces artificial adjustments that deviate from the fluid ``optimal'' path. To minimize such deviation from expected paths, the threshold should be smaller if the waiting time to the next re-solving epoch is large in relative terms, as later re-solving poses a larger risk of excessive deviation. Since the proportion of time left to the next re-solving epoch $(t_{k+1}^*-t_k^*)/(T-t_k^*)$ decreases over $k=0,\ldots,K-1$, thresholds increase accordingly, leading to more aggressive thresholding toward the end.

\subsection{Technical novelty}
Before ending this section, we would like to comment on some key distinctions between our problem and analyses from those in~\citet{BumpensantiWa20}, whose work in network revenue management inspired the $\birt$ policy. 

The most salient feature of our problem is that its objective is non-linear and hence non-separable over time. This renders inapplicable a direct decomposition of the time horizon into several, disjoint re-solving epochs and plugging in single-epoch regret guarantees (assuming memorylessness of the allocations). Rather, in our problem, each re-solving epoch starts with the random utilities garnered from stochastic arrivals in previous epochs, and these initial conditions have to be taken into account in epoch-wise analyses.

In addition, in the dynamic fair allocation problem, the central planner needs to assign resources across multiple agents, i.e., allocation decisions lie in the simplex $\Delta_n$ across the $n$ agents. In the network revenue management problem, the decision maker needs to decide whether to fulfill demand requests, i.e., allocation decisions are binary accept-reject. This discrepancy complicates the thresholding procedure as thresholding one agent's assignment affects the decisions of all other agents. 

In analyzing the performance, \citet{BumpensantiWa20} invoke existing results by~\citet{ReimanWa08} for network revenue management problems, whereas we provide a new and simpler stand-alone analysis that exploits the structure of our problem. Our analysis also facilitates a more fundamental understanding of the impact of stochastic arrivals on the optimal solution through a new sensitivity analysis and of the coupling between the hindsight optimal solution and auxiliary policies. 

\section{Experiments}
\label{sec:experiment}
In this section, we present numerical results from our simulation experiments to illustrate the performances of the fluid and $\birt$ policies. The first experiment evaluates these policies on several randomized arrival distributions under the egalitarian welfare, the harmonic-mean welfare and the Nash social welfare. The second experiment focuses on the special case of the egalitarian welfare and exposes the impact of degeneracy on re-solving heuristics. It is built upon two artificially designed arrival distributions, whose fluid problems are degenerate and nondegenerate respectively.

\subsection{Randomized experiment}
In the first experiment, we illustrate the performances of the $\F$ and $\birt$ policies under the egalitarian welfare, the harmonic-mean welfare and the Nash social welfare, using a number of randomized problem instances. More specifically, we consider a problem setup of $n=4$ agents and $L=5$ types of resources. We simulate probabilities of the resource types uniformly at random from the simplex, i.e., $(p_\ell:\ell\in[L])\sim\mathsf{Unif}(\Delta_L)$, and the marginal utilities independently and identically from a parameterized $\mathsf{Beta}$ distribution, i.e., $\beta_\ell^i \overset{\text{iid}}{\sim}\mathsf{Beta}(\cdot, \cdot)$ across $\ell, i$. In particular, we make the two choices for the marginal utilities:
\begin{itemize}
\item $\beta_\ell^i \sim \mathsf{Beta}(0.5, 0.5)$, under which values are concentrated around the two extremes of $0$ and $1$, or in other words, the agents have dichotomous valuations for the resources;
\item $\beta_\ell^i \sim \mathsf{Beta}(2.0, 2.0)$, under which values are concentrated near the middle at $0.5$, or in other words, the agents have similar valuations for the resources.
\end{itemize}
In each case, we simulate 30 instances of arrival distributions and compute the average relative regrets over the randomized arrival distributions. Note that relative regrets are used to normalize the problems for a reasonable comparison regardless of welfare dimensions. Results are shown in Figure~\ref{fig:randomized-experiments}. The meta-distributions for generating arrival distributions are distinguished by marks, and the policies are distinguished by colors.

{
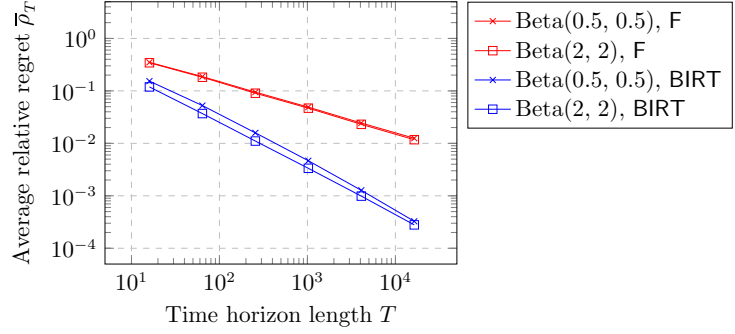
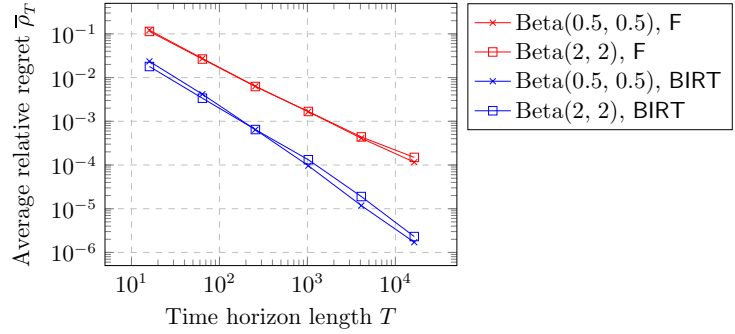
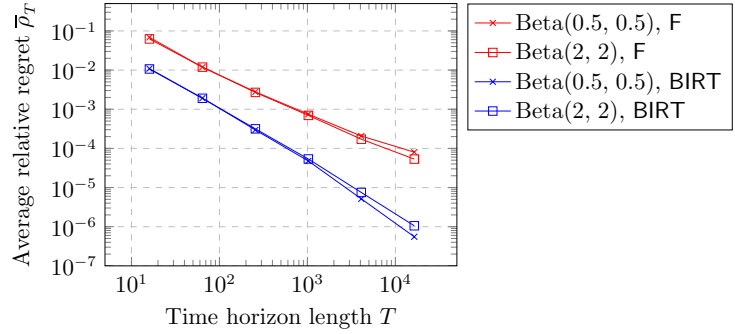
\begin{figure}[htp]
\centering
\begin{subfigure}[b]{0.9\textwidth}
\centering
\scalebox{0.8}{
\begin{tikzpicture}
\begin{loglogaxis}
[
    width=0.5\textwidth,
    height=0.4\textwidth,
    ylabel={Average relative regret $\overline\rho_T$},
    xlabel={Time horizon length $T$},
    xmin=5, xmax=50000,
    ymin=0.00005, ymax=5,
    xtick={10,100,1000,10000},
    ytick={0.0001, 0.001, 0.01, 0.1, 1, 10},
    legend pos=outer north east,
    legend cell align={left},
    xmajorgrids=true,
    ymajorgrids=true,
    grid style=dashed,
]

\addplot[color=red,mark=x,]
coordinates {
(16,	0.34997071256308265)
(64,	0.18735603739722717)
(256,	0.09465295779672023)
(1024,	0.049048659793681226)
(4096,	0.024396199618357914)
(16384,	0.01250270185334418)
};\addlegendentry{Beta(0.5, 0.5), $\mathsf{F}$}
\addplot[color=red,mark=square]
coordinates {
(16,	0.3437342345238327)
(64,	0.18204399285888637)
(256,	0.09029334347670195)
(1024,	0.04684914511939584)
(4096,	0.023045263409568015)
(16384,	0.01173382049273235)
};\addlegendentry{Beta(2, 2), $\mathsf{F}$}

\addplot[color=blue, mark=x,]
coordinates {
(16, 	 0.1535601864553649)
(64, 	 0.05250535317470937)
(256, 	 0.015858964738274018)
(1024, 	 0.00466424957894454)
(4096, 	 0.0012716567544070222)
(16384, 	 0.00032653715306897384)
};\addlegendentry{Beta(0.5, 0.5), $\mathsf{BIRT}$}
\addplot[color=blue,mark=square,]
coordinates {
(16, 	 0.1189992266434671)
(64, 	 0.036955392767107936)
(256, 	 0.011054781603075324)
(1024, 	 0.0033433951837953386)
(4096, 	 0.000986338686706168)
(16384, 	 0.00027893156573165234)
};\addlegendentry{Beta(2, 2), $\mathsf{BIRT}$}
\end{loglogaxis}
\end{tikzpicture}
}
\vspace{-10pt}
\caption{The egalitarian welfare ($q=-\infty$)}
\label{fig:randomized-experiment-egal}
\end{subfigure}
\\
\vspace{10pt}
\begin{subfigure}[b]{0.9\textwidth}
\centering
\scalebox{0.8}{
\begin{tikzpicture}
\begin{loglogaxis}
[
    width=0.5\textwidth,
    height=0.4\textwidth,
    ylabel={Average relative regret $\overline\rho_T$},
    xlabel={Time horizon length $T$},
    xmin=5, xmax=50000,
    ymin=5e-7, ymax=0.5,
    xtick={10,100,1000,10000},
    ytick={1e-7, 1e-6, 1e-5, 1e-4, 1e-3, 1e-2, 1e-1, 1},
    legend pos=outer north east,
    legend cell align={left},
    xmajorgrids=true,
    ymajorgrids=true,
    grid style=dashed,
]

\addplot[color=red,mark=x,] coordinates {
(16,	0.12117877405239247)
(64,	0.02735437367920858)
(256,	0.006325200394097331)
(1024,	0.001658903627589495)
(4096,	0.0004094301837426477)
(16384,	0.00011629368907769061)
};\addlegendentry{Beta(0.5, 0.5), $\mathsf{F}$}
\addplot[color=red,mark=square,]coordinates {
(16,	0.11381849386352022)
(64,	0.02655284785729576)
(256,	0.006239495278306232)
(1024,	0.0016762885878102568)
(4096,	0.00043901107324404505)
(16384,	0.00014935811358876973)
};\addlegendentry{Beta(2, 2), $\mathsf{F}$}

\addplot[color=blue, mark=x,] coordinates {
(16, 	 0.023497552123135382)
(64, 	 0.00416501548590546)
(256, 	 0.0006436997968224125)
(1024, 	 9.638770920551082e-05)
(4096, 	 1.181450557188559e-05)
(16384, 	 1.713509713573606e-06)
};\addlegendentry{Beta(0.5, 0.5), $\mathsf{BIRT}$}
\addplot[color=blue,mark=square,]coordinates {
(16, 	 0.017923338043610293)
(64, 	 0.003377021704319229)
(256, 	 0.0006473810570042528)
(1024, 	 0.00013188310279161418)
(4096, 	 1.8997823143333546e-05)
(16384, 	 2.3080633473229445e-06)
};\addlegendentry{Beta(2, 2), $\mathsf{BIRT}$}
\end{loglogaxis}
\end{tikzpicture}
}
\vspace{-10pt}
\caption{The harmonic-mean welfare ($q=-1$)}
\label{fig:randomized-experiment-harmonic}
\end{subfigure}
\\
\vspace{10pt}
\begin{subfigure}[b]{0.9\textwidth}
\centering
\scalebox{0.8}{
\begin{tikzpicture}
\begin{loglogaxis}
[
    width=0.5\textwidth,
    height=0.4\textwidth,
    ylabel={Average relative regret $\overline\rho_T$},
    xlabel={Time horizon length $T$},
    xmin=5, xmax=50000,
    ymin=1e-7, ymax=5e-1,
    xtick={10,100,1000,10000},
    ytick={1e-7, 1e-6, 1e-5, 1e-4, 1e-3, 1e-2, 1e-1, 1},
    legend pos=outer north east,
    legend cell align={left},
    xmajorgrids=true,
    ymajorgrids=true,
    grid style=dashed,
]

\addplot[color=red,mark=x,]coordinates {
(16,	0.06868510095688479)
(64,	0.011892379110217469)
(256,	0.0027705373865693638)
(1024,	0.0007472175126065768)
(4096,	0.00020673828321773714)
(16384,	8.011983154023501e-05)
};\addlegendentry{Beta(0.5, 0.5), $\mathsf{F}$}
\addplot[color=red,mark=square,]coordinates {
(16,	0.06251875947320219)
(64,	0.011904375538923019)
(256,	0.0026933989349730014)
(1024,	0.0006974945084220463)
(4096,	0.00017179952543922095)
(16384,	5.331174464874114e-05)
};\addlegendentry{Beta(2, 2), $\mathsf{F}$}

\addplot[color=blue, mark=x,]coordinates {
(16, 	 0.010892784389423742)
(64, 	 0.0019210666893304047)
(256, 	 0.0002958315156290427)
(1024, 	 4.767311723278842e-05)
(4096, 	 5.149232258565453e-06)
(16384, 	 5.530262106530742e-07)
};\addlegendentry{Beta(0.5, 0.5), $\mathsf{BIRT}$}
\addplot[color=blue,mark=square,]coordinates {
(16, 	 0.010690196782383882)
(64, 	 0.0019041282664236376)
(256, 	 0.00031576417226560045)
(1024, 	 5.36577881102756e-05)
(4096, 	 7.492319711388827e-06)
(16384, 	 1.0573789130130156e-06)
};\addlegendentry{Beta(2, 2), $\mathsf{BIRT}$}
\end{loglogaxis}
\end{tikzpicture}
}
\vspace{-10pt}
\caption{Nash social welfare ($q=0$)}
\label{fig:randomized-experiment-nash}
\end{subfigure}
\caption{Average relative regret $\overline\rho(\pi)$ over randomized arrival distributions, under the egalitarian welfare ($q=-\infty$), the harmonic-mean welfare ($q=-1$) and the Nash social welfare ($q=0$), plotted against the time horizon length $T$ on a double logarithmic scale}
\label{fig:randomized-experiments}
\end{figure}
}

We summarize our findings under the three welfare metrics.
\begin{itemize}
\item Figure~\ref{fig:randomized-experiment-egal} shows that under the egalitarian welfare ($q=-\infty$), the $\F$ policy attains a relative regret of $\rho_T(\F) = \Theta(T^{-1/2})$, translating to a regret of $\R_T(\F) = \Theta(\sqrt{T})$; the $\birt$ policy attains a relative regret of $\rho_T(\birt) = \Theta(T^{-1})$, translating to a uniformly bounded regret of  $\R_T(\birt) = \Theta(1)$.
\item Figure~\ref{fig:randomized-experiment-harmonic} shows that under the harmonic-mean welfare ($q=-1$), the $\F$ policy attains a relative regret of $\rho_T(\F) = \Theta(T^{-1})$, translating to a regret of $\R_T(\F) = \Theta(1)$; the $\birt$ policy attains a relative regret of $\rho_T(\birt) = o(T^{-1})$, translating to a vanishing regret of $\R_T(\birt) = o(1)$.
\item Figure~\ref{fig:randomized-experiment-nash} shows that under the Nash social welfare ($q=0$), the $\F$ policy attains a relative regret of $\rho_T(\F) = \Theta(T^{-1})$, translating to a regret of $\R_T(\F) = \Theta(1)$; the $\birt$ policy attains a relative regret of $\rho_T(\birt) = o(T^{-1})$, translating to a vanishing regret of $\R_T(\birt) = o(1)$.
\end{itemize}

These empirical findings from the randomized experiments exactly match out theoretical guarantees. Interestingly, we notice the $\birt$ policy attains a resounding $o(1)$ vanishing regret under the harmonic-mean and Nash social welfare metrics, an even more impressive improvement upon the $\Theta(1)$ regret attained by the static $\F$ policy. The precise mechanism to attaining the vanishing regret is unknown to us, but the empirical findings do highlight the excellent performance that $\birt$ has under all welfare metrics.

\subsection{Special experiment under the egalitarian welfare}
\label{sec:experiment-special}
In the second experiment, we zoom in and focus on the egalitarian welfare. We aim to corroborate our claims on the performances of the fluid and $\birt$ policies under the egalitarian welfare and, moreover, to understand what driving forces are behind the outstanding performance of $\birt$. To be more precise, we are including the following heuristic policies for comparison.
\begin{enumerate}
\item $\F$: the fluid static policy based on the initial fluid problem (Algorithm~\ref{alg:fluid}).
\item $\birt$: the Backward Infrequent Re-solving with Thresholding policies (Algorithm~\ref{alg:birt}).
\item $\mathsf{BIR}$: the Backward Infrequent Re-solving policies (Algorithm~\ref{alg:bir} in Appendix~\ref{sec:add-alg}), which re-solve the fluid problems on the same schedule as $\birt$ does but without thresholding (i.e., the updated fluid policy $\boldsymbol\xi^\F$ at re-solving epoch $t_k^*$ is directly used as $\boldsymbol\xi^k$ for implementation during the current epoch).
\item $\mathsf{FR}$: the Frequent Re-solving policy (Algorithm~\ref{alg:fr} in Appendix~\ref{sec:add-alg}), which re-solves the fluid problems in every period before making allocation decisions.
\end{enumerate}
We include the $\mathsf{BIR}$ policy (Algorithm~\ref{alg:bir}) to illustrate the crucial role that thresholding is playing in improving performance and the $\mathsf{FR}$ policy (Algorithm~\ref{alg:fr}) to evaluate whether performance could be improved by substituting thresholding with more frequent re-solving.

In addition to inter-policy comparison, we would like to investigate empirically what impact a degenerate deterministic problem has on the performance of re-solving heuristics, as widely discussed in previous literature. 

To illustrate the impact of degeneracy, we consider a simple problem setup with $n=2$ agents and $L=2$ types of resources where the marginal utilities are $\beta_1 = (1,\nicefrac12)$ and $\beta_2 = (\nicefrac12,1)$. We then consider two arrival distributions $(p_1,p_2) = (\nicefrac12, \nicefrac12)$ and $(p_1,p_2) = (\nicefrac{2}{5}, \nicefrac{3}{5})$. The first case (detailed in Appendix~\ref{sec:proof-opt-flu}) results in a degenerate fluid problem, as a total of $6>5$ constraints are active at the fluid policy $\boldsymbol\xi^\F$ with $(\xi^\F)_1=(1,0)$ and $(\xi^\F)_2=(0,1)$, corresponding to the two bases $(u, \xi_1^1, \xi_2^1, \xi_2^2)$ and $(u, \xi_1^1, \xi_1^2, \xi_2^2)$. Therefore, slight perturbations on the number of arrivals can lead to optimal hindsight allocations of the form $(\xi)_1=(1,0)$ and $(\xi)_2=(\epsilon,1-\epsilon)$ or $(\xi)_1=(1-\epsilon,\epsilon)$ and $(\xi)_2=(0,1)$ for $\epsilon \ge 0$, and, without thresholding, the resolving policy might allocate resources to the wrong agent (e.g., allocating resource one to agent two when the hindsight optimal basis is the first one). In the second case, the original fluid problem is nondegenerate as a total of $5$ constraints are active at the fluid policy $\boldsymbol\xi^\F$ with $(\xi^\F)_1=(1,0)$ and $(\xi^\F)_2=(\nicefrac{2}{9},\nicefrac{7}{9})$, corresponding to the unique basis $(u, \xi_1^1, \xi_2^1, \xi_2^2)$. Therefore, with high probability, an optimal hindsight allocation is of the form $(\xi)_1=(1,0)$ and $(\xi)_2=(\nicefrac{2}{9}\pm\epsilon,\nicefrac{7}{9}\mp\epsilon)$ and it is less likely to allocate resources to the wrong agent. In both cases, the same hyper-parameter $\eta=1.05$ is chosen for the $\birt$ and $\mathsf{BIR}$ policies. Regrets of respective policies are shown in Figure~\ref{fig:regret}.

{
\begin{figure}[htbp]
\centering
\begin{minipage}[t]{0.48\textwidth}
\centering
\begin{tikzpicture}[scale=0.9]
\begin{loglogaxis}[
    width=\textwidth,
    ylabel={Regret $\mathcal{R}_T(\pi)$},
    xlabel={Time horizon length $T$},
    xmin=10, xmax=100000,
    ymin=0.05, ymax=100,
    xtick={10,100,1000,10000, 100000},
    ytick={0.1, 1, 10, 100},
    axis y line=left,
    legend pos=outer north east,
    legend cell align={left},
    xmajorgrids=true,
    ymajorgrids=true,
    grid style=dashed,
]

\addplot[
    color=red,
    mark=*,
    ]
    coordinates {
(16, 	 1.0413333330341017)
(64, 	 2.282666665663766)
(256, 	 4.324666663757095)
(1024, 	 8.641999970549275)
(4096, 	 16.3819998358397)
(16384, 	 33.801332606105504)
(65536, 	 69.64533220725713)
    };
    \addlegendentry{$\mathsf{F}$}

\addplot[
    color=orange,
    mark=triangle,
    ]
    coordinates {
(16, 	 0.25467555502786565)
(64, 	 0.46174888725469676)
(256, 	 1.0092866630129285)
(1024, 	 2.28625997323323)
(4096, 	 4.60910653009596)
(16384, 	 8.897850308829245)
(65536, 	 21.85312357058768)
    };
    \addlegendentry{$\mathsf{FR}$}
    
\addplot[
    color=blue,
    mark=square,
    ]
    coordinates {
(16, 	 0.28703133264759484)
(64, 	 0.49473366512189265)
(256, 	 1.0394871622147641)
(1024, 	 1.8384044678340943)
(4096, 	 4.084565827080258)
(16384, 	 7.291233102280944)
(65536, 	 14.969055699118366)
    };
    \addlegendentry{$\mathsf{BIR}$}

\addplot[
    color=black,
    mark=o,
    ]
    coordinates {
(16, 	 0.2836748326488812)
(64, 	 0.43443316515341196)
(256, 	 0.6842921623229905)
(1024, 	 0.8635769695771405)
(4096, 	 1.1221463354416372)
(16384, 	 1.1677671294403016)
(65536, 	 1.1103162417773529)
    };
\addlegendentry{$\mathsf{BIRT}$}

\end{loglogaxis}
\end{tikzpicture}%
\end{minipage}
\hfill
\begin{minipage}[t]{0.48\textwidth} 
\centering
\begin{tikzpicture}[scale=0.9]
\begin{loglogaxis}[
    width=\textwidth,
    xlabel={Time horizon length $T$},
    xmin=10, xmax=100000,
    ymin=0.05, ymax=100,
    xtick={10,100,1000,10000, 100000},
    ytick={0.1, 1, 10, 100},
    axis y line=right,
    xmajorgrids=true,
    ymajorgrids=true,
    grid style=dashed,
]

\addplot[
    color=red,
    mark=*,
    ]
    coordinates {
(16, 	 1.175335999190654)
(64, 	 2.669758663877716)
(256, 	 5.156098328130565)
(1024, 	 10.370699658679023)
(4096, 	 19.862195971764002)
(16384, 	 40.34064522090073)
(65536, 	 83.03905655051142)
    };

\addplot[
    color=orange,
    mark=triangle,
    ]
    coordinates {
(16, 	 0.3687555546926595)
(64, 	 0.6276044419572007)
(256, 	 0.5952633294529599)
(1024, 	 0.633128881662371)
(4096, 	 0.8390544160522662)
(16384, 	 0.7141798880975997)
(65536, 	 0.7064306605018403)
    };
    
\addplot[
    color=blue,
    mark=square,
    ]
    coordinates {
(16, 	 0.3883034991906582)
(64, 	 0.6142081638777273)
(256, 	 0.7396968281305631)
(1024, 	 0.8122061586790039)
(4096, 	 0.8665819717639015)
(16384, 	 0.8746107209002238)
(65536, 	 0.8110705505113438)
    };

\addplot[
    color=black,
    mark=o,
    ]
    coordinates {
(16, 	 0.38226049919065824)
(64, 	 0.5826831638777279)
(256, 	 0.6842798281321345)
(1024, 	 0.7283611586790041)
(4096, 	 0.7684644717639021)
(16384, 	 0.7651497209002563)
(65536, 	 0.730853050511414)
    };

\end{loglogaxis}
\end{tikzpicture}%
\end{minipage}

\caption{Comparison of regrets of $\mathsf{F}$, $\mathsf{BIR}$, $\mathsf{FR}$, and $\mathsf{BIRT}$ policies on a double logarithmic scale. Results for the degenerate instance are shown on the left, and those for the nondegenerate instance are shown on the right.}
\label{fig:regret}
\end{figure}
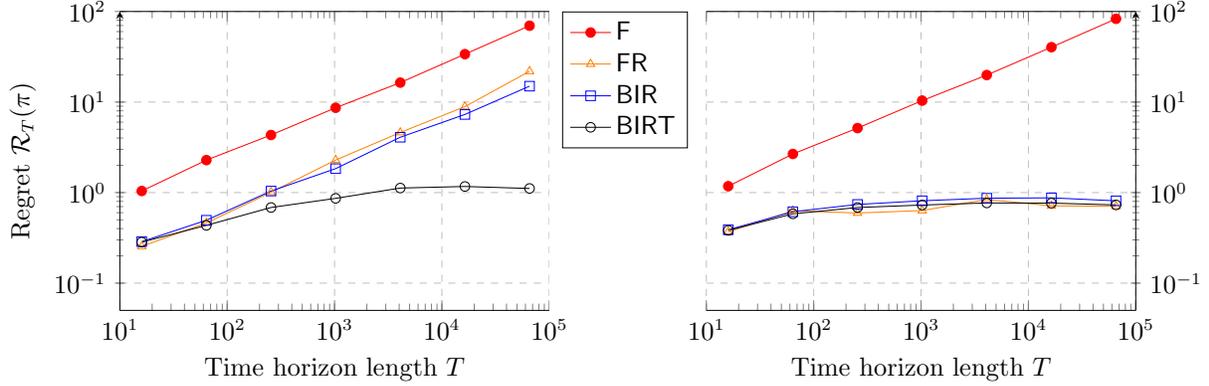
}
First, we observe from Figure~\ref{fig:regret} that in both cases $\birt$'s regret converges, corroborating our claim (Theorem~\ref{thm:main}) that it attains $O(1)$ regret. Equipped with distributional knowledge, the central planner is shown to be able to improve the performance upon the fluid static policy by re-solving infrequently ($\mathsf{BIR}$ policy), either by improving the factor of $\sqrt{T}$ in the regret in the degenerate case, or by attaining a constant regret in the nondegenerate case. Interestingly, in both cases, further increasing the re-solving frequency does not lead to performance improvement; in fact, the $\mathsf{BIR}$ and $\mathsf{FR}$ policies attain similar regrets. In other words, carefully designing the re-solving schedule suffices to guarantee good performance while drastically reducing computational complexity. The most significant difference between the two cases is that the $\mathsf{BIR}$ and $\mathsf{FR}$ policies attain $\Theta(\sqrt{T})$ regret in the degenerate case, whereas they can attain $O(1)$ regret in the nondegenerate case, indicating their vulnerability against degeneracy. By contrast, the more robust $\birt$ policy steadily attains a $O(1)$ regret in both cases, demonstrating that the thresholding rule is crucial for controlling the excessive variation in allocations and obtain remarkable performance even in the nondegenerate case.

\section{Conclusion}
\label{sec:conclusion}

In this work we study dynamic policies for allocating $T$ divisible items arriving sequentially to a fixed set of $n$ agents. 
We consider a collection of H\"older-mean welfare metrics, including the Nash social welfare and the egalitarian welfare, and evaluate policies based on the expected welfare generated, which entails addressing the trade-off between efficiency and fairness. We show the classical fluid policy based on certainty equivalents incurs a $\Theta(\sqrt{T})$ regret against the expected hindsight optimal egalitarian welfare and an $O(1)$ regret under all other welfare metrics. For the egalitarian welfare, we propose the Backward Infrequent Re-solving with Thresholding ($\birt$) policies, which incur an $O(1)$ regret, uniformly bounded over time horizon length $T$. 

It would be interesting to extend the work to problems defined on an infinite time horizon, or with infinitely many arrival types. To allow for flexibility in the arrival model, an open question is to design dynamic allocation policies that are robust against inexact solutions to the fluid problems and distribution mis-specification in the IID case (inaccurate utilities), and that work well when arrivals are inter-temporally dependent (e.g., autocorrelated, evolving under a Markov chain model).

\bibliographystyle{informs2014} 
\bibliography{bib} 

\ECSwitch


\ECHead{Appendices}

\section{Proofs}
\label{sec:proofs}

\subsection{Proof of Lemma~\ref{lemma:opt<=flu}}
\label{sec:proof-opt<=flu}
Fix $P\in\P$. By Lemma~\ref{lemma:opt-static}, $P$-almost surely, there is a static policy $\boldsymbol\xi^*$ that is hindsight optimal, i.e.,$\opt = w (\sum_{\ell\in[L]} N_\ell \beta_\ell * \xi^*_\ell )$. Consider an alternate static policy $\boldsymbol\xi \in \Delta_n^L$ given by $\xi_\ell^i := \EE[N_\ell (\xi^*)_\ell^i] / \EE[N_\ell]$, for all $i, \ell$. Hence,
\begin{equation*}
\begin{split}
    \EE[\opt]  = 
    \EE\left[ w \left(\sum_{\ell\in[L]} N_\ell \beta_\ell * \xi^*_\ell \right) \right]  \le w \left(\sum_{\ell\in[L]} \beta_\ell^i\EE[N_\ell (\xi^*)_\ell^i]  \right) 
     = w \left(\sum_{\ell\in[L]}  \EE[N_\ell]\beta_\ell^i \xi_\ell^i \right)
    \le \flu,
\end{split}
\end{equation*}
where the first equality follows by definition of $\boldsymbol\xi^*$, the first inequality by Jensen's inequality, the second equality by definition of $\boldsymbol\xi$, and the second inequality by definition of $\flu$.

\subsection{Proof of Lemma~\ref{lemma:flu-flu}}
\label{sec:proof-flu-flu}
Fix static policy $\boldsymbol\xi$ and welfare metric $w\equiv w_q$---we omit the notation for convenience in the proof. Denote $\underline{p}:=\min_\ell p_\ell$. Clearly, we can assume $w(\EE[B_T]) > 0$ because the lemma is otherwise trivial.
First notice 
\[w(B_T) \ge w\left(\EE[B_T] \min_i\q{B_T^i}{\EE[B_T^i]} \right) = w(\EE[B_T]) \min_i 
\q{B_T^i}{\EE[B_T^i]} , \] where the inequality follows because $w$ is nonnegative and monotone, and the equality because $w$ is positively homogeneous. Hence,  
\[1-\q{\EE[w(B_T)]}{w(\EE[B_T])} \le \EE\left[\max_\ell \q{\EE[B_T^i]-B_T^i}{\EE[B_T^i]}\right] \le \sqrt{\EE\left(\max_\ell \q{\EE[B_T^i]-B_T^i}{\EE[B_T^i]} \right)^2} \le \sqrt{\sum_\ell \EE \left(\q{\EE[B_T^i]-B_T^i}{\EE[B_T^i]}\right)^2} \le \sqrt{\q{L}{\underline{p}T}},\]
where the first inequality follows because $w(\EE[B_T])>0$, the second by Jensen's inequality, the third by nonnegativity of squares, and the last by Lemma~\ref{lemma:variance}. Since $w$ is positively homogeneous, $w(\EE[B_T]) = T w(\sum_{\ell\in[L]} p_\ell\beta_\ell*\xi_\ell)$, so \[w(\EE[B_T]) - \EE[w(B_T)] \le \sqrt{\q{LT}{\underline{p}}} w \left( \sum_{\ell\in[L]} p_\ell\beta_\ell*\xi_\ell \right) \le \q1n \sqrt{\q{LT}{\underline{p}}},\]
where the last inequality follows by Lemma~\ref{lemma:as-bound}. The proof is thereby concluded.
\Halmos

\subsection{Proof of Theorem~\ref{thm:main-fluid}}
\label{sec:proof-main-fluid}
Fix static policy $\boldsymbol{\xi}$ and omit the notation for convenience. We show the following claim, to which the theorem is a corollary by invoking Lemma~\ref{lemma:opt<=flu}.

\textbf{Claim}: Under any H\"older-mean welfare metric $w_q$ with $q\in(-\infty,1]$ and any static policy $\boldsymbol\xi$, we have $w_q(\EE[B_T]) - \EE[w_q(B_T)] = O(1)$ uniformly bounded over time horizon length $T$.

We assume $w_q(\EE[B_T])>0$ because the theorem is trivial otherwise. For convenience, we denote $\underline{p}:=\min_\ell p_\ell >0$ and define $Z^i := B_T^i / \EE [B_T^i] - 1$ for each agent $i$ if $\EE[B_T^i]>0$ and $Z^i=0$ otherwise. Then $\EE [Z^i] = 0$ and $Z^i\in[-1,1/\underline p]$ a.s.\ because $B_T^i\ge0$ and
\begin{equation}\label{eqn:z-as-bound}
1+Z^i = \q{B_T^i}{\EE[B_T^i]} = \q{\sum_\ell N_\ell\beta_\ell^i \xi_\ell^i}{\sum_\ell T p_\ell \beta_\ell^i \xi_\ell^i} \le \q{\sum_\ell T \beta_\ell^i \xi_\ell^i}{\sum_\ell T p_\ell \beta_\ell^i \xi_\ell^i} \le \q1{\underline{p}} \quad\text{ a.s.}
\end{equation}
because $N_\ell \le T$ for all $\ell \in [L]$ and $p_\ell \ge \underline p$.

In the sequel, we show the claim separately for $q\in(0,1)$, $q=0$ and $q\in(-\infty,0)$.

{\proof{Proof of Theorem~\ref{thm:main-fluid} for $q\in(0,1)$.}
First notice Lemma~\ref{lemma:taylor-power-lower-bound} implies
\[\q{w_q(B_T)}{w_q(\EE[B_T])} 
\ge 1+\q1q \q{\sum_{i\in[n]} (\EE[B_T^i])^{q}\left[(1+Z^i)^{q}-1\right]}{\sum_{i\in[n]} (\EE[B_T^i])^{q}} 
\ge 1+\q{\sum_{i\in[n]} (\EE[B_T^i])^{q}\left[Z^i - (q^{-1} - 1)(Z^i)^2\right]}{\sum_{i\in[n]} (\EE[B_T^i])^{q}},\] 
so taking expectations yields
\[1 - \q{\EE[w_q(B_T)]}{w_q(\EE[B_T])} 
\le (q^{-1}-1)\q{\sum_{i\in[n]} (\EE[B_T^i])^{q} \EE[(Z^i)^2]}{\sum_{i\in[n]} (\EE[B_T^i])^{q}} 
\le \q{q^{-1}-1}{\underline{p}}T^{-1},\]
where the first inequality follows because $\EE[Z^i]=0$, and the second by Lemma~\ref{lemma:variance}. \endproof}

{\proof{Proof of Theorem~\ref{thm:main-fluid} for Nash social welfare, i.e., $q=0$.}
Consider the event $E_{1/2}=\{|Z^i|\le 1/2, \forall i\}$, on which
\[ \q{w_q(B_T)}{w_q(\EE[B_T])} \ge 1 + \log\q{w_q(B_T)}{w_q(\EE[B_T])} = 1 + \q1n\sum_{i\in[n]} \log(1+Z^i) \ge 1 + \q1n \sum_{i\in[n]} \left[Z^i - (Z^i)^2\right],\]
where the first inequality follows because $\exp(x)\ge1+x$ for all $x\in\RR$, the equality by definition of $w$ and $Z^i$, and the second inequality by Lemma~\ref{lemma:taylor-log}. We then notice 
\[\begin{split}
\EE\left[\q{w_q(B_T)}{w_q(\EE[B_T])}\right] 
&\ge \EE\left[\q{w_q(B_T)}{w_q(\EE[B_T])}; E_{1/2}\right] \\
&\ge \EE\left[1 + \q1n \sum_{i\in[n]} Z^i - \q1n\sum_{i\in[n]} \left(Z^i\right)^2; E_{1/2} \right] \\
&\ge \EE \left[1 + \q1n \sum_{i\in[n]} Z^i; E_{1/2} \right] - \q1n\sum_{i\in[n]} \EE\left(Z^i\right)^2 \\
&= 1 - \EE \left[\q1n \sum_{i\in[n]} (1+Z^i); E_{1/2}^\complement \right] - \q1n\sum_{i\in[n]} \EE\left(Z^i\right)^2,
\end{split}\]
where the first inequality follows by the almost sure nonnegativity of $w_q(B_T)$, the second inequality follows by the almost sure nonnegativity of squares, and the equality because $\EE[Z^i]=0$. Hence, \begin{equation*}
1 - \EE\left[\q{w_q(B_T)}{w_q(\EE[B_T])}\right] \le \EE \left[\q1n \sum_i (1+Z^i); E_{1/2}^\complement \right] + \q1n\sum_i\EE(Z^i)^2 \le \q{2L}{\underline{p}} \exp\left(-\q12\underline{p}^2T\right) + (\underline{p}T)^{-1}
\end{equation*} because Equation~\eqref{eqn:z-as-bound} implies $1+Z^i \le p^{-1}$ a.s., Lemma~\ref{lemma:hoeffding} implies $\PP(E_{1/2}^\complement) \le 2L\exp(-\underline{p}^2T/2)$, and Lemma~\ref{lemma:variance} implies $\EE(Z^i)^2 \le 1/\underline{p}T$. \endproof}

{\proof{Proof of Theorem~\ref{thm:main-fluid} for $q\in(-\infty,0)$.} On the event $E_{1/2}=\{|Z^i|\le 1/2, \forall i\}$, $B_T>0$, and 
\[\q{w_q(B_T)}{w_q(\EE[B_T])} \ge 1+\q1q \q{\sum_{i\in[n]} (\EE[B_T^i])^{q}\left[(1+Z^i)^{q}-1\right]}{\sum_{i\in[n]} (\EE[B_T^i])^{q}} \ge 1+\q{\sum_{i\in[n]} (\EE[B_T^i])^{q}\left[Z^i + 2^{-q+2}(Z^i)^2/q\right]}{\sum_{i\in[n]} (\EE[B_T^i])^{q}},\]
where the first inequality follows by Lemma~\ref{lemma:taylor-power-lower-bound}, and the second by Lemma~\ref{lemma:taylor-power-upper-bound}. We then notice 
\[\begin{split}
\EE\left[\q{w_q(B_T)}{w_q(\EE[B_T])}\right] 
&\ge \EE\left[\q{w_q(B_T)}{w_q(\EE[B_T])}; E_{1/2}\right] \\
&\ge \EE\left[1 + \q{\sum_{i\in[n]} (\EE[B_T^i])^{q}\left[Z^i + 2^{-q+2}(Z^i)^2/q\right]}{\sum_{i\in[n]} (\EE[B_T^i])^{q}} ; E_{1/2} \right] \\
&\ge \EE \left[\q{\sum_{i\in[n]} (\EE[B_T^i])^{q}(1+Z^i)}{\sum_{i\in[n]} (\EE[B_T^i])^{q}}; E_{1/2} \right] + \q{2^{-q+2}}{q} \q{\sum_{i\in[n]} (\EE[B_T^i])^{q} \EE(Z^i)^2}{\sum_{i\in[n]} (\EE[B_T^i])^{q}}\\
&= 1 - \EE \left[\q{\sum_{i\in[n]} (\EE[B_T^i])^{q}(1+Z^i)}{\sum_{i\in[n]} (\EE[B_T^i])^{q}}; E_{1/2}^\complement \right] + \q{2^{-q+2}}{q} \q{\sum_{i\in[n]} (\EE[B_T^i])^{q} \EE(Z^i)^2}{\sum_{i\in[n]} (\EE[B_T^i])^{q}}\\
&\ge 1 - \EE\left[\max_i(1+Z^i); E_{1/2}^\complement\right] + \q{2^{-q+2}}{q} \max_i \EE(Z^i)^2,
\end{split}\]
where the first inequality follows by the almost sure nonnegativity of $w_q(B_T)$, the second inequality by above, the third by the almost sure nonnegativity of squares, the equality because $\EE[Z^i]=0$, and the last inequality because the coefficients add up to one. Hence, \begin{equation*}
1 - \EE\left[\q{w_q(B)}{w_q(\EE[B_T])}\right] \le \EE \left[\max_i (1+Z^i); E_{1/2}^\complement \right] + \q{2^{-q+2}}{-q} \max_i\EE(Z^i)^2 \le \q{2L}{\underline{p}} \exp\left(-\q12\underline{p}^2T\right) + \q{2^{-q+2}}{(-q)\underline{p}}T^{-1}
\end{equation*} because Equation~\eqref{eqn:z-as-bound} implies $1+Z^i \le p^{-1}$ a.s., Lemma~\ref{lemma:hoeffding} implies $\PP(E_{1/2}^\complement) \le 2L\exp(-\underline{p}^2T/2)$, and Lemma~\ref{lemma:variance} implies $\EE(Z^i)^2 \le 1/\underline{p}T$. \endproof}

Finally, we conclude by arguing $w_q(\EE[B_T]) = Tw_q(\sum_\ell p_\ell\beta_\ell*\xi_\ell)$ because $w_q$ is positively homogeneous for any $q\in(-\infty,1)$ and by invoking Lemma~\ref{lemma:opt<=flu}. \Halmos

\subsection{Proof of Lemma~\ref{lemma:opt-flu}}
\label{sec:proof-opt-flu}
Consider the arrival distribution $P$ where $P\{b=(1,\nicefrac12)\} = \nicefrac12 = P\{b=(\nicefrac12,1)\}$. The hindsight optimum is $\opt = \nicefrac T2 - |N_1-N_2|/6$. WLOG suppose $N_1\le N_2$, then this is attained by the offline policy $\xi_1=(1,0)$ and $\xi_2=(4N_2-2T, 2T-N_2)/3N_2$. The fluid policy is $(\xi^\F)_1=(1,0)$ and $(\xi^\F)_2=(0,1)$, which achieves only $\alg(\F) = \nicefrac T2 - |N_1-N_2|/2$. Hence, Lemma~\ref{lemma:1/2} implies
\[\R_T(\F) = \EE[\opt - \alg(\F)] = \q23 \EE\left|N_1-\q T2\right| = \Theta(\sqrt{T}).\]

\subsection{Proof of Theorem~\ref{thm:main}}
\label{sec:proof-main}
We start with the telescopic decomposition of the regret:
\begin{equation}
\label{eqn:telescope}
    \R_T(\birt) = \EE\left[\opt-\alg(\birt)\right] 
     = \sum_{k=1}^{K}\EE\left[\opt^{k-1}-\opt^k\right] + \EE\left[\opt^K - \alg(\birt)\right],
\end{equation}
where we recall $\opt=\opt^0$, and $\opt^k$ is the hindsight optimal egalitarian welfare starting with welfare $B_{t_k^*}$ at time $t_k^*$. First notice the last term is subject to an almost sure uniform bound
\begin{equation}
\label{eqn:uniform2}
\opt^K - \alg(\birt) \le T-t_K^* = 2, 
\end{equation}
implied by Lemma~\ref{lemma:uniform} as if the last re-solving occurs at $T$, because $\opt^K$ and $\birt$ agree on all actions by $t_K^* = T-2$. Now we proceed with the first terms. For $k=1,\ldots,K$, the policies $\opt^{k-1}$ and $\opt^k$ agree on all allocations by $t_{k-1}^*$, so Proposition~\ref{prop:one-epoch} implies (as if time is relabeled to starting at $t_{k-1}^*$)
\[ \EE\left[\opt^{k-1} - \opt^k | B_{t_{k-1}^*} \right] \le 3L(T-t_{k-1}^*) \exp\left(-C\q{(T-t_{k}^*)^4}{(T-t_{k-1}^*)^3} \right).\]
Notice the constant $C$ stays the same for all $k$, as it is determined only by the distribution $P$ and is independent of time $T$ and initial welfare $B_{t_{k-1}^*}$. To simplify the right hand side, we lower bound $T-t_{k}^*$ and upper bound $T-t_{k-1}^*$ as follows.

First, $T - t_{k}^* = \lfloor \exp(\eta^{K-k}) \rfloor \ge \lfloor \mathrm{e} \rfloor = 2$ implies \[\q{T-t_{k}^*}{\exp(\eta^{K-k})} = \q{\lfloor \exp(\eta^{K-k}) \rfloor}{\exp(\eta^{K-k})} \ge \q{\lfloor \exp(\eta^{K-k}) \rfloor}{\lfloor \exp(\eta^{K-k}) \rfloor + 1} \ge \q23.\]
Then notice $T-t_{k-1}^* = \lfloor \exp(\eta^{K-k+1}) \rfloor \le \exp(\eta^{K-k+1})$ for $k=2,\ldots,K$. For the corner case $k=1$, $T-t_0^* = T \le \exp(\eta^{K})$ by definition of $K = \lceil\log\log T/\log \eta \rceil$. Hence, denoting $C':=(\nicefrac23)^4C$, we have for any $k=1,\ldots, K$,
\[\EE\left[\opt^{k-1} - \opt^k \right] \le 3L \mathrm{e}^{\eta^{K-k+1}} \exp\left(-C' \mathrm{e}^{(4-3\eta)\eta^{K-k}} \right) .\]
Then we know
\begin{equation*}
\begin{split}
    \R_T(\birt) &\le 2 + 3L \sum_{k=1}^{K} \mathrm{e}^{\eta^{K-k+1}}  \exp\left(-C' \mathrm{e}^{(4-3\eta)\eta^{K-k}} \right) \\
    & = 2+3L\sum_{k=1}^{K} \mathrm{e}^{\eta^{k}} \exp\left(-C' \mathrm{e}^{(4-3\eta)\eta^{k-1}} \right) \\
    & \le 2+3L\sum_{k=1}^{\infty} \mathrm{e}^{\eta^{k}} \exp\left(-C' \mathrm{e}^{(4\eta^{-1}-3)\eta^{k}} \right),
\end{split}
\end{equation*}
where the first inequality follows from the telescopic decomposition~\eqref{eqn:telescope} and the almost sure uniform bound~\eqref{eqn:uniform2}, the equality by re-ordering the terms, and the last inequality by non-negativity.

Finally, we invoke Lemma~\ref{lemma:series}, a result on the sum of the infinite series, to conclude that
\[\R_T(\birt) = \EE[\opt-\alg(\birt)] \le 2 +  3LC'^{-\eta/(4-3\eta)} h\left(\eta,4\eta^{-1}-3\right) = O(1).\]

\subsection{Proof of Proposition~\ref{prop:one-epoch}}
\label{sec:proof-one-epoch}
We first give the following uniform bound on the difference between $\opt$ and $\opt^1$, before proceeding to bound the probability that the difference is nonzero. The proof is deferred to Appendix~\ref{sec:proof-uniform}.
\begin{lemma}
\label{lemma:uniform}
Suppose all arrival utilities are bounded $b_t\in[0, 1]$ for all $t$, and $\opt^1$ is the hindsight optimal egalitarian welfare starting with utilities $B_{t_1}$ at time $ t_1\in[T]$, with $B_{t_1} \ge B_0$. Then its egalitarian welfare loss against the hindsight optimum starting with utilities $B_0$ at time $0$ is at most
\begin{equation*}
    \opt - \opt^1 \le t_1.
\end{equation*}
\end{lemma}

Now it remains to bound the decoupling probability $\PP\{\opt\ne \opt^1\}$, because
\[ \EE\left[\opt - \opt^1|B_0\right] = \EE\left[(\opt-\opt^1)\textbf{1}\{\opt\ne \opt^1\}|B_0\right] \le t_1 \PP\{\opt\ne \opt^1|B_0\}. \]

We next provide a general sufficient condition on the coupling event $\{\opt=\opt^1\}$. For ease of notation, we denote the number of type-$\ell$ arrivals by $t_1$ by $N_\ell^\le := \sum_{t=1}^{t_1} \textbf{1}\{b_t = \beta_\ell\}$, and that of those after $t_1$ by $N_\ell^> := \sum_{t=t_1+1}^T \textbf{1}\{b_t = \beta_\ell\}$. We denote using tilted notation the deviation from their respective means (e.g., $\tilde{N}_\ell := N_\ell - \EE[N_\ell]$). The proof is deferred to Appendix~\ref{sec:proof-con}.
\begin{lemma}
\label{lemma:con}
Suppose the online player acts according to a static policy $\boldsymbol\xi^0\in\Delta_n^L$ for $t=1,2,\ldots, t_1$. Then $\opt=\opt^1$ if there is a hindsight optimal static solution $\boldsymbol\xi^*\in\Delta_n^L$ such that
\begin{equation}
\label{eqn:con}
N_\ell (\xi^*)_\ell \ge N_\ell^\le (\xi^0)_\ell, \quad\forall \ell.
\end{equation}
\end{lemma}

Now that we have a general sufficient condition for $\opt=\opt^1$, the remaining parts of the proof can be outlined as the following three steps: 1) invoking sensitivity analysis results for linear programming to show the existence of a hindsight optimal policy that is close to the fluid policy, 2) providing sufficient conditions on the arrivals under which the problem $\opt$ and $\opt^1$ are equivalent, and 3) proving that these sufficient conditions hold with high probability using Hoeffding’s inequality.

We start by showing that there exists a hindsight optimal policy that is close to the fluid policy, by invoking the LP sensitivity analysis literature against right-hand side constants and exploiting the similarity between the offline problem and the fluid problem. Since the $\birt$ policy is constructed based on the fluid policy, this guarantees the $\birt$ policy is close to the hindsight optimal policy, i.e., stochastic fluctuations are small compared to the leeway that the prophet $\opt^1$ has after time $t_1$ to make up for previous suboptimal allocations. Hence, we will be able to show that $\opt^1$ can recover the exact same utilities for all agents as $\opt$ does with high probability.

We first give a lemma on the proximity between fluid policies and hindsight optimal policies. The proof is deferred to Appendix~\ref{sec:proof-sensitivity2}.
\begin{lemma}
\label{lemma:sensitivity2}
For any fluid policy $\xi^\F$, there exists an hindsight optimal policy $\xi^*$ such that 
$$\max_{i,\ell}\left| N_\ell (\xi^*)_\ell^i - Tp_\ell(\xi^\F)_\ell^i \right| \le (nL+1)S \norm{\tilde N}_\infty,$$
where $S\ge1$ is a constant determined solely by $(\beta_\ell)_{\ell\in[L]}$, the \textit{support} of the arrival distribution $P$.
\end{lemma}

Now we detail, in Lemma~\ref{lemma:con2}, sufficient conditions specific to thresholded fluid policies under which Equation~\eqref{eqn:con} is satisfied, so that Lemma~\ref{lemma:con} can be used to argue $\opt=\opt^1$. The proof can be found in Appendix~\ref{sec:proof-con2}.
\begin{lemma}
\label{lemma:con2}
Suppose $\xi^\F$ is a fluid policy, and the online player acts $\xi^0$ constructed by applying the thresholding rule with $\gamma \in [0, (T-t_1)/2n^2T]$, for $t=t_1+1, \ldots, T$. Then $\opt = \opt^1$ if 
\begin{equation}
\label{eqn:con2}
    \tilde N_\ell^\le \le \q{(T-t_1)p_\ell}2 \text{ for all }\ell,\quad\text{ and }\quad \norm{\tilde N}_\infty \le \q{(T-t_1)\underline p}2 \q{\gamma}{S(nL+1)},
\end{equation}
where we define $\underline p:= \min_\ell p_\ell$.
\end{lemma}

Using the union bound, Lemma~\ref{lemma:con2} implies 
\[ \PP\left\{\opt\ne\opt^1|B_0\right\} \le \sum_{\ell\in[L]} \PP\left\{\tilde N_\ell^\le > \q{(T-t_1)p_\ell}2\right\} + \sum_{\ell\in[L]} \PP\left\{\left|\tilde N_\ell\right| > \q{(T-t_1)\underline p}2 \q{\gamma}{S(nL+1)}\right\}. \]

We then proceed to invoke Hoeffding's inequality to bound the probabilities. With $N_\ell^\le \sim \mathrm{Bin}(t_1, p_\ell)$ and $N_\ell \sim \mathrm{Bin}(T, p_\ell)$,  we have for all $\ell\in[L]$,
\[ \PP\left\{\tilde N_\ell^\le > \q{(T-t_1)p_\ell}2\right\} \le \exp\left(-\q{(T-t_1)^2}{2t_1}p_\ell^2 \right), \quad \text{ and} \]
\[ \PP\left\{\left|\tilde N_\ell\right| > \q{(T-t_1)\underline p}2 \q{\gamma}{S(nL+1)}\right\} \le 2 \exp\left(-\q{(T-t_1)^2}{2T} \q{\gamma^2}{S^2(nL+1)^2} \underline p^2 \right). \]
Hence,
\begin{equation*}
    \PP\left\{\opt\ne\opt^1|B_0\right\} \le L \left[ \exp\left(-\q{(T-t_1)^2}{2t_1}\underline p^2 \right) + 2\exp\left(-\q{(T-t_1)^2}{2T} \q{\gamma^2}{S^2(nL+1)^2} \underline p^2 \right) \right].
\end{equation*}

In particular, we can take the threshold $\gamma = (T-t_1)/2n^2T$ for the tightest bound, so that
\begin{equation*}
    \PP\left\{\opt\ne\opt^1|B_0\right\} \le L \left[ \exp\left(-\q{\underline p^2}2\q{(T-t_1)^2}{t_1} \right) + 2\exp\left(-\q{\underline p^2}{8n^4S^2(nL+1)^2}\q{(T-t_1)^4}{T^3} \right) \right].
\end{equation*}
Since $(T-t_1)^2/t_1 \ge (T-t_1)^4/T^3$ and $S\ge1$, we have
\begin{equation*}
    \PP\left\{\opt\ne\opt^1|B_0\right\} \le 3L \exp\left(-\q{\underline p^2}{8n^4S^2(nL+1)^2}\q{(T-t_1)^4}{T^3} \right),
\end{equation*}
and the proof is concluded upon invoking Lemma~\ref{lemma:uniform}.

\subsection{Proof of Lemma~\ref{lemma:uniform}}
\label{sec:proof-uniform}

Suppose a hindsight optimal policy allocates $\boldsymbol x^* \in \Delta_n^T$. Then
\begin{equation*}
\begin{split}
    \opt - \opt^1 
    &= \min_{1\le j \le n} \left( B_0^j + \sum_{t=1}^T b_t^j (x^*)_t^j \right) - \max_{(x_t)_t} \min_i \left( B_{ t_1}^i + \sum_{t= t_1 + 1}^T b_t^i x_t^i \right) \\
    &\le \min_{1\le j \le n} \left( B_0^j + \sum_{t=1}^T b_t^j (x^*)_t^j \right) - \min_{i\in[n]} \left( B_{ t_1}^i + \sum_{t= t_1 + 1}^T b_t^i (x^*)_t^i \right) \\
    &= \max_{i\in[n]} \left[ \min_{1\le j \le n} \left( B_0^j + \sum_{t=1}^T b_t^j (x^*)_t^j \right) - \left( B_{ t_1}^i + \sum_{t= t_1 + 1}^T b_t^i (x^*)_t^i \right) \right] \\
    &\le \max_{i\in[n]} \left( B_0^i - B_{ t_1}^i +\sum_{t=1}^{ t_1} b_t^i (x^*)_t^i \right) \\
    &\le \max_{i\in[n]} \sum_{t=1}^{ t_1} b_t^i (x^*)_t^i \\
    &\le t_1,
\end{split}
\end{equation*}
where the first equality follows by definition of $\boldsymbol x^*$, the first inequality because $\boldsymbol x^* \in \Delta_n^T$, the second equality because $-\min = \max$, the second inequality by taking the minimizing agent $i$ in the latter term, the third inequality because $B_0 \le B_{ t_1}$, and the last inequality because $b_t^i\le1$ and $(x^*)_t^i \le 1$.

\subsection{Proof of Lemma~\ref{lemma:con}}
\label{sec:proof-con}

Clearly, $\opt=\opt^1$ if for any arrival type, any agent can get as much utility from $\opt^1$'s policy as from $\opt$'s policy, i.e.\ there exists $\boldsymbol\xi\in\Delta_n^L$ such that 
$$N_\ell (\xi^*)_\ell^i = N_\ell^\le (\xi^0)_\ell^i + N_\ell^> \xi_\ell^i.$$

Fix arrival type $\ell$. If $N^>_\ell=0$, then we claim Equation~\eqref{eqn:con} holds with equality. Otherwise, $N_\ell = \sum_i N_\ell (\xi^*)_\ell^i > \sum_i N_\ell^\le (\xi^0)_\ell^i = N_\ell^\le = N_\ell$, which is absurd. Now we have shown $N_\ell (\xi^*)_\ell^i = N_\ell^\le (\xi^0)_\ell^i$ for all $i$, i.e., all agents have accumulated the exact welfare needed from type-$\ell$ arrivals by time $t_1$. An arbitrary allocation for type $\ell$ would satisfy the condition, since there will be no such arrivals.
    
If $N^>_\ell >0$, the proof is concluded by simply taking the static policy $\boldsymbol\xi\in\Delta_n^L$ given by
$$\xi_\ell = \q{N_\ell (\xi^*)_\ell - N_\ell^\le (\xi^0)_\ell}{N_\ell^>}, \quad \forall \ell\in[L].$$

\subsection{Proof of Lemma~\ref{lemma:sensitivity2}}
\label{sec:proof-sensitivity2}

We restate the following result on sensitivity of optimal solutions against RHS in linear programs.
\begin{lemma}[\citet{Schrijver98}, Theorem 10.5] 
\label{lemma:sensitivity}
Let $A$ be an $p\times m$-matrix, and let $S$ be such that for each nonsingular submatrix $M$ of $A$ all entries of $M^{-1}$ are at most $S$ in absolute value. Let $c$ be a row $n$-vector, and let $b'$ and $b''$ be column $m$-vectors such that $\max \{cx|Ax\le b'\}$ and $\max \{ c x | A x \le b''\}$ are finite. Then for each optimum solution $x'$ of the first maximum there exists an optimum solution $x''$ of the second maximum with
$$\norm{x'-x''}_\infty \le mS \norm{b'-b''}_\infty.$$
\end{lemma}

To apply the sensitivity analysis result, we rewrite in the linear program form the offline problem as (with optimal solution $(y^*)_\ell^i = N_\ell (\xi^*)_\ell^i$)
\begin{align*}
    \opt = \max \quad U &&& \\
    \text{s.t.}\quad U &\le B_0^i + \sum_\ell \beta_\ell^i y_\ell^i, & & \forall i, \\
    y_\ell^i & \ge0, && \forall i,\ell,\\
    \sum_i y_\ell^i & = N_\ell, && \forall \ell,
\end{align*}
and the initial fluid problem can be rewritten as (with optimal solution $(z^*)_\ell^i := T p_\ell (\xi^\F)_\ell^i$)
\begin{align*}
    \flu = \max \quad U &&& \\
    \text{s.t.}\quad U &\le B_0^i + \sum_\ell \beta_\ell^i z_\ell^i, & & \forall i, \\
    z_\ell^i & \ge0, && \forall i,\ell,\\
    \sum_i z_\ell^i & = T p_\ell, && \forall \ell.
\end{align*}
Lemma~\ref{lemma:sensitivity} implies for any optimal solution $\mathbf z$ to the fluid problem, there exists an optimal solution $\mathbf y$ to the offline problem such that 
$$\max\left\{\left|\flu-\opt\right|, \norm{\mathbf{y-z}}_\infty\right\} \le (nL+1)S \norm{\tilde N}_\infty,$$
where $S$ is larger than any absolute value of entries in the inverse of any nonsingular submatrix of 
\begin{equation*}
A = \left[
\begin{array}{l|c}
    \textbf{1}_n & \begin{array}{ccc}
        \mathsf{diag}(-\beta_1) & \cdots & \mathsf{diag}(-\beta_L)\end{array} \\ \hline
    \textbf{0}_{n\times L} & -I_{n\times L} \\ \hline
    \textbf{0}_L & \begin{array}{ccc}
        \textbf{1}_n^\top & \cdots & 0 \\
        \vdots & \ddots & \vdots \\
        0 & \cdots & \textbf{1}_n^\top 
    \end{array} \\ \hline
    \textbf{0}_L & \begin{array}{ccc}
        -\textbf{1}_n^\top & \cdots & 0 \\
        \vdots & \ddots & \vdots \\
        0 & \cdots & -\textbf{1}_n^\top 
    \end{array}
\end{array}
\right].
\end{equation*}
Here the decision variable $x$ in the statement of the lemma is $x = (U, y_1, \ldots, y_L)$ or $x = (U, z_1, \ldots, z_L)$ and the right-hand side vector of the constraints is $(B_0, \textbf{0}_{n\times L}, \vec N, -\vec N)$ or $(B_0, \textbf{0}_{n\times L}, T \vec p, -T \vec p)$, respectively.  In particular, we also observe $S\ge1$. The proof is concluded as a corollary.

\subsection{Proof of Lemma~\ref{lemma:con2}}
\label{sec:proof-con2}

First, Lemma~\ref{lemma:sensitivity2} implies the existence of a hindsight optimal static policy $\boldsymbol\xi^* \in \Delta_n^L$ such that 
$$\max_{i,\ell}\left| N_\ell (\xi^*)_\ell^i - Tp_\ell(\xi^\F)_\ell^i \right| \le (nL+1)S \norm{\tilde N}_\infty,$$
so for all $i$ and $\ell$,
\begin{equation}
\label{eqn:diff}
\begin{split}
    N_\ell (\xi^*)_\ell^i - N_\ell^\le (\xi^\F)_\ell^i 
    &\ge \left[ Tp_\ell (\xi^\F)_\ell^i - (nL+1)S \norm{\tilde N}_\infty\right] - \left[ t_1 p_\ell + \tilde N_\ell^\le \right] (\xi^\F)_\ell^i \\
    &= \left[ (T-t_1) p_\ell - \tilde N_\ell^\le \right] (\xi^\F)_\ell^i - (nL+1)S \norm{\tilde N}_\infty.
\end{split}
\end{equation}
Lemma~\ref{lemma:con} implies that it suffices to show $N_\ell(\xi^*)_\ell^i \ge N_\ell^\le (\xi^0)_\ell^i$ for all $i$ and $\ell$. Now we check condition~\eqref{eqn:con} element-wise. Recall from the thresholding rule that $(\xi^0)_\ell^i \in \{0\} \cup[\gamma,1]$ for all $i$ and $\ell$. 
\begin{itemize}
\item 
If $(\xi^0)_\ell^i=0$, then $N_\ell (\xi^*)_\ell^i \ge 0 = N_\ell^\le (\xi^0)_\ell^i$.
    
\item
If $(\xi^0)_\ell^i\ge\gamma$ and $i\ne j$, then $(\xi^\F)_\ell^i = (\xi^0)_\ell^i \ge \gamma$. We have
\begin{equation*}
\begin{split}
N_\ell (\xi^*)_\ell^i - N_\ell^\le (\xi^0)_\ell^i &= N_\ell (\xi^*)_\ell^i - N_\ell^\le (\xi^\F)_\ell^i \\
&\ge \left[ (T-t_1) p_\ell - \tilde N_\ell^\le \right] (\xi^\F)_\ell^i - (nL+1)S \norm{\tilde N}_\infty \\
&\ge \q{(T-t_1) p_\ell}2 (\xi^\F)_\ell^i - \q{(T-t_1) p_\ell}2 \gamma\\
&\ge0,
\end{split}
\end{equation*}
where the equality follows because $(\xi^\F)_\ell^i = (\xi^0)_\ell^i$, the first inequality by Equation~\eqref{eqn:diff}, the second by assumption~\eqref{eqn:con2}, and the last because $(\xi^\F)_\ell^i = (\xi^0)_\ell^i \ge \gamma$.

\item
Agent $j\in \arg\max_i (\xi^\F)_\ell^i$ is the only one receiving positive adjustment. We have
\begin{equation*}
\begin{split}
    N_\ell (\xi^*)_\ell^j - N_\ell^\le (\xi^0)_\ell^j &= N_\ell (\xi^*)_\ell^j - N_\ell^\le (\xi^\F)_\ell^j - N_\ell^\le (\xi^0 - \xi^\F)_\ell^j \\
    &\ge N_\ell (\xi^*)_\ell^j - N_\ell^\le (\xi^\F)_\ell^j - n\gamma N_\ell^\le \\
    &\ge \left[ Tp_\ell (\xi^\F)_\ell^j - (nL+1)S\norm{\tilde N}_\infty\right] - N_\ell^\le (\xi^\F)_\ell^j - n\gamma N_\ell^\le \\
    &= \left[ (T-t_1)p_\ell - \tilde N_\ell^\le \right]\left[ (\xi^\F)_\ell^j + n\gamma \right] - (nL+1) S \norm{\tilde N}_\infty - n\gamma T p_\ell \\
    &\ge \q{(T-t_1)p_\ell}2 \left[ (\xi^\F)_\ell^j + n\gamma \right] - \q{(T-t_1)p_\ell}2\gamma - n\gamma T p_\ell \\
    &\ge nTp_\ell \left[ \q{T-t_1}{2n^2T} - \gamma \right] + \q{(T-t_1)p_\ell}2 (n-1)\gamma\\
    &\ge0.
\end{split}
\end{equation*}
where the first equality follows by rearranging terms, the first inequality because the adjustment is at most $(\xi^0-\xi^\F)_\ell^j = -\sum_{i\ne j} (\xi^0-\xi^\F)_\ell^i < n\gamma$, the second inequality by Equation~\eqref{eqn:diff}, the second equality because $N_\ell^\le = t_1p_\ell + \tilde N_\ell^\le$, the third inequality by assumption~\eqref{eqn:con2}, the fourth inequality because $(\xi^\F)_\ell^j = \max_{i}(\xi^\F)_\ell^i \ge \nicefrac1n$, and the last inequality by the threshold rule $0\le \gamma_0 \le (T-t_1)/2n^2T$.
\end{itemize}
Therefore, we conclude $N_\ell(\xi^*)_\ell^i \ge N_\ell^\le (\xi^0)_\ell^i$ for all $i$ and $\ell$, and so $\opt=\opt^1$ by Lemma~\ref{lemma:con}.

\section{Additional results}
\label{sec:add-results}

\begin{lemma}[Futility of randomization]
\label{lemma:futility}
Consider the online convex optimization problem
\begin{equation*}
\begin{split}
    \mathrm{maximize}\quad & w\left(\sum_{t=1}^T b_t x_t \right)\\
    \mathrm{subject\,to}\quad & x_t \in \X_t,\\
    & x_t \perp \sigma(b_{t+1},\ldots,b_T)
\end{split}
\end{equation*}
where $w:\RR^m \to \RR$ is a concave function, $b_t\in\RR^{m\times n}$ are fixed arrivals, and $\X_t\subseteq \RR^n$ is a convex set for every $t$. Then in expectation, any randomized online policy is weakly dominated by some deterministic online policy.
\end{lemma}
{\proof{Proof of Lemma~\ref{lemma:futility}}
Fix any randomized online policy $\pi$, and suppose it outputs random actions $(x_t(\pi):t\in[T])$. Jensen's inequality implies that
\begin{equation*}
     \EE_{\pi} \left[w\left(\sum_{t=1}^T b_t x_t(\pi) \right)\right] \le w\left(\sum_{t=1}^T b_t \EE_{\pi} \left[x_t(\pi)\right] \right).
\end{equation*}
Notice the mean policy $\bar\pi$ given by actions $x_t(\bar\pi) = \EE_{\pi} \left[x_t(\pi)\right]$ is feasible, because $x_t(\bar\pi)\in\X_t$ by the convexity of the action domains $\X_t$, so it is an online policy, and the RHS is simply its objective.
\Halmos\endproof}

\begin{lemma}[Static policies dominate]
\label{lemma:opt-static}
$P$-almost surely, any offline policy is weakly dominated by a static policy.
\end{lemma}
{\proof{Proof of Lemma~\ref{lemma:opt-static}}
Based on any offline policy $\boldsymbol{x}$, the proof is immediate by defining the static policy given by 
\begin{equation*}
    \xi_\ell := \q1{N_\ell}\sum_{t:b_t=\beta_\ell} x_t \quad\text{ for all } \ell\in[L]\text{ where }N_\ell>0,
\end{equation*}
because this static policy satisfies $\xi_\ell \in\Delta_n$ and gives agents the same utilities as the offline policy does.
\Halmos\endproof}

\begin{lemma}
\label{lemma:prophet-fair}
$P$-almost surely, there is a hindsight optimal policy under which the agents have equal utilities. 
\end{lemma}
{\proof{Proof of Lemma~\ref{lemma:prophet-fair}}
Fix an arrival sequence $\boldsymbol{b}\equiv(b_t: t\in[T])$. We first deal with the corner case $\opt=0$. This means that a na"ive even allocation ($x_t^i=1/n$ for all $t$ and $i$) also attains zero welfare, so there is some agent $i^*$ where $b_t^{i^*}=0$ for all $t$. Hence, allocating all resources to agent $i^*$ results in all agents receiving zero utility.

In the sequel suppose the hindsight optimum $\opt>0$ is achieved by some optimal offline allocation sequence $\x$ that result in unequal utilities across agents. Then we can denote the set of agents with excess utilities as $\mathcal{I}:=\{i: \sum_t b_t^i x_t^i > \opt\} \ne \varnothing$. Notice by definition, $\sum_t b_t^i x_t^i \ge \opt$ for all $i\in[n]$, so $\sum_t b_t^i x_t^i = \opt$ for $i\in[n]\setminus \mathcal{I}$. We aim to construct an alternative allocation sequence that equalize all agents' utilities while attaining the hindsight optimum $\opt$.

We first consider the heuristic allocation sequence $\hat\x$ given by 
\[\hat x_t^{i} := \q{\opt}{\sum_t b_t^i x_t^i} x_t^{i} + \q1n \sum_{j\in[n]} \left(1-\q{\opt}{\sum_t b_t^j x_t^j}\right)x_t^{j},\quad\quad i=1,\ldots,n.\]
Clearly $\hat x_t \in\Delta_n$, because $\hat x_t\ge\mathbf0$ and $\sum_i \hat x_t^i = 1$. Then under $\hat\x$, agent $i$ would have a cumulative utility of 
\[\sum_{t=1}^T b_t^{i}\hat x_t^{i} 
= \opt + \q1n \sum_{j\in[n]} \left(1-\q{\opt}{\sum_t b_t^j x_t^j}\right) \sum_{t=1}^T b_t^{i} x_t^{j}.\]
By definition, $\hat\x$ cannot perform better than the optimal offline policy, so $\min_{1\le i\le n} \sum_{t=1}^T b_t^{i}\hat x_t^{i} \le \opt$. Since $\sum_t b_t^i x_t^i \ge \opt$ for all $i$ with strict inequality holding for $j\in\mathcal{I}$, this implies there exists some agent ${i^*} \in [n]$ such that
\begin{equation}
\label{eqn:cs}
    b_t^{({i^*})} x_t^{j} = 0, \quad \forall j\in\mathcal{I}, t\in[T].
\end{equation}

This necessary condition~\eqref{eqn:cs} is similar to the well-known complementary slackness condition in duality theory. It fixes an agent $i^*$ and stipulates that resources with positive utility to $i^*$ ($t$ when $b_t^{i^*}>0$) may not be allocated to agents with excess utilities, i.e., $x_t^{j}=0$ for all $j\in\mathcal{I}$. In particular, we notice that $i^*\not\in\mathcal{I}$, since otherwise $b_t^{i^*} x_t^{i^*}=0$ for all $t$, contradicting $\sum_t b_t^{i^*} x_t^{i^*} \ge \opt>0$. This means $i^*$ is a binding agent, i.e., $\sum_t b_t^{i^*} x_t^{i^*} = \min_{1\le i\le n}\sum_t b_t^i x_t^i = \opt$.

After we pick out the agent ${i^*}$, condition~\eqref{eqn:cs} implies that any positive allocation $x_t^{j}>0$ to agents with excess utilities $j\in\mathcal{I}$ can be safely redirected toward agent $i^*$ without increasing its cumulative utility. 

More formally, consider the allocation sequence $\tilde\x$ given by
\begin{equation*}
\begin{split}
    \tilde x_t^{i} &:= \q{\opt}{\sum_t b_t^i x_t^i} x_t^{i},\quad\quad\quad \forall i\ne{i^*},\\
    \tilde x_t^{i^*} &:=  x_t^{i^*} +  \sum_{j\in[n]}  \left(1-\q{\opt}{\sum_t b_t^j x_t^j}\right) x_t^{j}.
\end{split}
\end{equation*}
Clearly $\tilde x_t \in\Delta_n$, because $\tilde x_t\ge\mathbf0$, and $\sum_i x_t^i = 1$. Furthermore, it is in fact a resource-equalizing hindsight optimal solution, because
\begin{equation*}
\begin{split}
    \sum_{t=1}^T b_t^{i}\tilde x_t^{i} 
    &= 
    \sum_{t=1}^T b_t^{i} \q{\opt}{\sum_t b_t^i x_t^i} x_t^i
    = 
    \opt, \quad\quad\forall i\ne i^*,\\
    \sum_{t=1}^T b_t^{i^*}\tilde x_t^{i^*} &= \sum_{t=1}^T b_t^{i^*}x_t^{i^*} +  \underbrace{\sum_{j\in\mathcal{I}}  \left(1-\q{\opt}{\sum_t b_t^j x_t^j}\right) \sum_{t=1}^T b_t^{i^*}x_t^{j}}_{= 0\text{ by definition of }i^*~\eqref{eqn:cs}} = \sum_t b_t^{i^*} x_t^{i^*} = \opt,
\end{split}
\end{equation*}
where the last inequality follows because $i^*\not\in\mathcal{I}$. Hence we have constructed a hingsight optimal allocation that equalize the agents' cumulative utilities.
\Halmos\endproof}

\begin{lemma}
\label{lemma:flu-opt}
There exists some $P\in\P$ under which $\flu - \EE[\opt] = \Omega(\sqrt{T})$.
\end{lemma}
{\proof{Proof of Lemma~\ref{lemma:flu-opt}}
Take a simple example where $n=2$ agents start with zero utilities $B_0=(0,0)$, and the arrival structure is given by $\beta_1 = (1,0)$ and $\beta_2 = (0,1)$ with respective probabilities $p_1=\nicefrac12$ and $p_2=\nicefrac12$. Clearly the optimal value of the fluid problem is $\flu = \nicefrac T2$, attained by the fluid policy $\xi_1 = (1,0)$ and $\xi_2 = (0,1)$. 

To compute $\opt$, note the trivial bound $\opt \le \min(\sum_{t=1}^T b_t^{1}, \sum_{t=1}^T b_t^{2}) = \min(N_1, N_2)$; in fact, it can be attained by the na\"ive policy $x_t = b_t$, so $\opt = \min(N_1, N_2)$. We then know
\begin{equation*}
    \flu - \EE\left[\opt\right] 
    = \EE\left| N_1 - \q T2 \right| = \Theta(\sqrt{T}),
\end{equation*}
by invoking Lemma~\ref{lemma:1/2} in Appendix~\ref{sec:add-proofs}.
\Halmos\endproof}

\begin{lemma}
\label{lemma:as-bound}
Under H\"older-mean welfare $w_q$ for $q\in[-\infty,1]$ and any static policy $\boldsymbol\xi$, 
$w_q(B_T) \le T / n$.
\end{lemma}
{\proof{Proof of Lemma~\ref{lemma:as-bound}} Under static policy $\boldsymbol\xi$,
\begin{equation*}
\begin{split}
w_q (B_T) 
& = w_q \left( \sum_{\ell\in[L]} N_\ell \beta_\ell * \xi_\ell \right) \\
& \le w_q \left( \sum_{\ell\in[L]} N_\ell \xi_\ell \right) \\
& \le w_q \left( \textbf{1} \q{\sum_{i\in[n]} \sum_{\ell\in[L]} N_\ell \xi_\ell^i}n \right) \\
& = \q{T}{n},
\end{split}
\end{equation*}
where the first equality follows by definition, the first inequality because $\beta_\ell \in [0,1]^n$, the second inequality because $w$ is concave, and the second equality because $\sum_\ell N_\ell=T$, $\xi_\ell\in\Delta_n$ and $w_q(\textbf{1})=1$.
}

\begin{lemma}\label{lemma:variance}
For any static policy $\xi$ and agent $i$, \[\Var(B_T^i) \le \q{\EE[B_T^i]^2}{\underline{p}T}.\]
\end{lemma}
{\proof{Proof of Lemma~\ref{lemma:variance}.}
First notice
\[\begin{split}
\Var(B_T^i) &= \sum_\ell \Var(N_\ell) (\beta_\ell^i \xi_\ell^i)^2 + \sum_{\ell\ne\ell'} \mathsf{Cov}(N_\ell, N_{\ell'}) \beta_\ell^i\beta_{\ell'}^i \xi_\ell^i \xi_{\ell'}^i \\
&= T \sum_\ell p_\ell(1-p_\ell) (\beta_\ell^i \xi_\ell^i)^2 + T \sum_{\ell\ne\ell'} (-p_\ell p_{\ell'}) \beta_\ell^i\beta_{\ell'}^i \xi_\ell^i \xi_{\ell'}^i \\
&= T \sum_\ell p_\ell (\beta_\ell^i \xi_\ell^i)^2 - T \left(\sum_{\ell} p_\ell \beta_\ell^i \xi_\ell^i \right)^2 \\
&\le T \sum_\ell p_\ell (\beta_\ell^i \xi_\ell^i)^2,
\end{split}\]
where the first equality follows because $B_T^i = \sum_\ell N_\ell \beta_\ell^i\xi_\ell^i$ by definition, the second equality because $(N_\ell)_\ell \sim \mathsf{Multinomial}(T, (p_\ell)_\ell)$, the third equality by rearranging terms, and the inequality by nonnegativity of squares. Next, note \[
\EE[B^i_T(\xi)]^2
= T^2 \left(\sum_\ell p_\ell\beta_\ell^i \xi_\ell^i \right)^2 
\ge T^2 \sum_\ell \left(p_\ell\beta_\ell^i \xi_\ell^i \right)^2 
\ge \underline{p} T^2 \sum_\ell p_\ell \left(\beta_\ell^i \xi_\ell^i \right)^2,\]
where the equality follows by definition, the first inequality follows by expanding the squares, and the second inequality because $p_\ell \ge \underline{p}$ for all $\ell$. The proof is thereby concluded.\Halmos
\endproof}

\begin{lemma}\label{lemma:hoeffding}
Given any $T\ge1$ and static policy $\xi$, define the event $E_\delta := \{|\tilde B_T^i| \le \delta \EE[B_T^i], \forall i\}$ for all $\delta>0$. Then $\PP(E_\delta^\complement) \le 2 L\exp(-2\underline{p}^2T\delta^2)$, where $\underline{p} := \min_\ell p_\ell$.
\end{lemma}
{\proof{Proof of Lemma~\ref{lemma:hoeffding}.}
\[\begin{split} 
\PP(E_\delta^\complement) 
&= \PP\left\{|B_T^i| > \delta \EE[B_T^i] \text{ for some } i \right\} \\
&\le \PP\left\{|\tilde N_\ell|>\delta \EE[N_\ell] \text{ for some } \ell \right\} \\
&\le \sum_{\ell\in[L]} \PP\{|\tilde N_\ell|>\delta \EE[N_\ell]\} \\
& \le \sum_{\ell\in[L]} 2\exp\left(-\q{2(\delta\EE[N_\ell])^2}{T}\right) \\
& = \sum_{\ell\in[L]} 2\exp\left(-2p_\ell^2T\delta^2\right) \\
&\le 2 L\exp(-2\underline{p}^2T\delta^2), 
\end{split}\]
where the first inequality follows because $|\tilde N_\ell| \le \delta \EE[N_\ell]$ for all $\ell$ implies $|B_T^i| \le \delta \EE[B_T^i]$ for all $i$, the second by the union bound, the third by Hoeffding's inequality, and the last because $p_\ell\ge\underline{p}$.
\Halmos\endproof}

We give a generalization of Theorem~\ref{thm:main} below, allowing for dependency of the distribution on time horizon length $T$.

\begin{definition}
\label{def:P-extension}
Given a finite support $\boldsymbol\beta \equiv(\beta_\ell \in [0,1]^n : \ell\in[L])$ and a positive $\underline{p} \in (0,\nicefrac12)$, the class of \emph{restricted admissible} probability distributions $\P(\underline p, \boldsymbol\beta)$ is the set of all finite probability distributions $P$ with support $\boldsymbol\beta$ such that $P(b = \beta_\ell) \ge \underline{p}$ for all $\ell\in[L]$.
\end{definition}

\begin{theorem}
\label{thm:extension}
(Extension of Theorem~\ref{thm:main}) Fix $\varepsilon\in(0,\nicefrac12)$. Fix any sequence $\underline p(T)\in(0,\nicefrac12)$ with $\underline p(T) = \Omega(T^{-\nicefrac12+\varepsilon})$. Then for any $\eta\in(1, (1-\nicefrac\varepsilon2)^{-1})$, 
\begin{equation}
    \sup_{P\in\P(\underline{p}(T), \boldsymbol\beta)} \R_T(\birt_\eta) = O(1).
\end{equation}
\end{theorem}
{\proof{Proof of Theorem~\ref{thm:extension}}
By assumption there exist $T_0>1$ and $\sigma > 0$ such that $\underline p(T) \ge \sigma T^{-\nicefrac12+\varepsilon}$ for all $T>T_0$. It suffices to assume $T>T_0$. Recall Proposition~\ref{prop:one-epoch} shows
\begin{equation*}
\begin{split}
    \PP\left\{\opt\ne\opt^1\right\} 
    &\le 3L \exp\left(-\q{\underline p^2}{8n^4S^2(nL+1)^2}\q{(T-t_1)^4}{T^3} \right) \\
    &\le 3L \exp\left(-\q{\sigma^2}{8n^4S^2(nL+1)^2}\q{(T-t_1)^4}{T^{4-2\varepsilon}} \right).
\end{split}
\end{equation*}
In particular, for re-solving points defined in the $\birt_\eta$ policy, we can show, similar to the proof of Theorem~\ref{thm:main} (Appendix~\ref{sec:proof-main}), that 
\[\EE\left[\opt^{k-1} - \opt^k \right] \le 3L \mathrm{e}^{\eta^{K-k+1}} \exp\left(-C'' \mathrm{e}^{(4-(4-2\varepsilon)\eta)\eta^{K-k}} \right),\]
for some constant $C''>0$ determined by $\sigma$, $T_0$, and the \textit{support} of the distribution. Finally, we invoke Lemma~\ref{lemma:series} to conclude, following a similar convergence argument as in the proof of Theorem~\ref{thm:main} (Appendix~\ref{sec:proof-main}), that for any $P\in\P(\underline{p}(T), \boldsymbol\beta)$ and any choice of $\eta\in(1,(1-\nicefrac\varepsilon2)^{-1})$,
\[\R_T(\birt) \le 2 + 3L C''^{-1/\zeta}h(\eta,\zeta) = O(1),\]
where $\zeta = 4\eta^{-1}-(4-2\varepsilon)\in(0,1)$.
\Halmos\endproof}

\section{Additional proofs}
\label{sec:add-proofs}
\begin{lemma}
\label{lemma:1/2}
Suppose $N_T \sim \mathsf{Bin}(T, \nicefrac{1}{2})$ for $T\ge1$. Then
\begin{equation*}
    \q{\EE|N_T-\nicefrac{T}{2}|}{\sqrt{T}} \to \q1{\sqrt{2\pi}} \quad\text{ as }\quad T\to\infty.
\end{equation*}
\end{lemma}
{\proof{Proof of Lemma~\ref{lemma:1/2}}
Define the sequence $M_T:=(2N_T-T)/\sqrt{T}$ for $T\ge1$. Then it suffices to show that $\EE|M_T|\to \sqrt{2/\pi}$ as $T\to\infty$.

Observe that $\EE M_T=0$ and $\EE M_T^2 = 4\mathrm{Var}(N_T)/T = 1$ for all $T$, so $\{M_T\}$ is uniformly integrable, and that the Central Limit Theorem implies $M_T\Rightarrow \mathsf{N}(0,1)$. By Skorokhod's Representation Theorem, there exists a sequence of random variables $(Z_T:T\ge1)$ such that $M_T\overset{d}{=}Z_T$ and $Z_T\to Z$ almost surely, where $Z\sim\mathsf{N}(0,1)$. This means the class $\{Z_T:T\ge1\}$ is also uniformly integrable, so $Z_T\to Z$ in $L^1$, implying in particular $\EE|Z_T| \to \EE|Z| = \sqrt{2/\pi}$. Hence, we conclude that $\EE|M_T| = \EE|Z_T| \to \sqrt{2/\pi}$.
\Halmos\endproof}

\begin{lemma}\label{lemma:taylor-log}
For any $x\ge-1/2$, $\log(1+x) \ge x - x^2$.
\end{lemma}
{\proof{Proof of Lemma~\ref{lemma:taylor-log}.}
Denote the smooth function $f \in \mathcal{C}^\infty([-1,\infty))$ given by $f(x) = \log(1+x) - x + x^2$, then $f(0)=0$ and $f(-1/2) = \nicefrac{3}{4}-\log2>0$. The derivative of $f(x)$ is \[
f'(x) = (1+x)^{-1} - 1 + 2x, 
\]
so $f'(0)=0$ and $f'(-1/2)=0$. Next, we argue that $0$ is the only local minimum of $f$.
The second-order derivative of $f$ is \[
f''(x) = -(1+x)^{-2} + 2.\]
Note it is monotonically increasing in $(-1, \infty)$; in particular, $f''<0$ on $(-1,-1+1/\sqrt{2})$ and $f''>0$ on $(-1+1/\sqrt{2},+\infty)$. Now $f'(x) = f'(0) + \int_0^x f''(z)\ud z$ for any $x\in(-1/2, \infty)$, so for any $x\ne0$ where $f'(x)=0$, $f''(x) < 0$. This means it cannot be a local minimum. Hence, $0$ is the only local minimum of $f$.

Finally, notice $f'(0)=0$ and $f''(x)>0$ for any $x\ge0$ implies that $f(x)>f(0)=0$ for any $x\ge0$. Since we know $f(0)=0$, we conclude that $0$ is a global minimum of $f$ and that $f(x)\ge0$ for all $x\ge-1/2$. \Halmos\endproof}

\begin{lemma}\label{lemma:taylor-power-lower-bound}
For $x\ge-1$, \[
(1+x)^q \ge \begin{cases}
1+qx - (1-q)x^2, & q\in(0,1),\\
1+qx, & \text{o/w}.
\end{cases}\]
\end{lemma}
{
\proof{Proof of Lemma~\ref{lemma:taylor-power-lower-bound}}
Suppose $q\in(0,1)$ and denote the smooth function $f \in \mathcal{C}^\infty([-1,\infty))$ given by $f(x) = (1+x)^{q} - 1 - qx + (1-q)x^2$, then $f(-1) = 0$ and $f(0)=0$. The derivative of $f(x)$ is \[
f'(x) = q (1+x)^{q-1} - q + 2(1-q) x, 
\]
so $f'(0)=0$ and $f'(-1^+)=+\infty$. Next, we argue that $0$ is the only local minimum of $f$.
The second-order derivative of $f$ is \[
f''(x) = q(q-1) (1+x)^{q-2} + 2(1-q).\]
Note it is monotonically increasing in $(-1, \infty)$, $f''(-1^+) = -\infty$, and $f''(0) = q^2-2q+2 > 0$. Now $f'(x) = f'(0) + \int_0^x f''(z)\ud z$ for any $x\in[-1, \infty)$, so for any $x\ne0$ where $f'(x)=0$, $f''(x) < 0$. This means it cannot be a local minimum. Hence, $0$ is the only local minimum of $f$.

Finally, notice $f'(0)=0$ and $f''(x)>0$ for any $x\ge0$ implies that $f(x)>f(0)=0$ for any $x\ge0$. Since we know $f(0)=0$, we conclude that $0$ is a global minimum of $f$ and that $f(x)\ge0$ for all $x\ge-1$. 

Now suppose $q\in(-\infty, 0]\cup[1,\infty)$ and consider the smooth function $g\in\mathcal{C}^\infty([-1,\infty))$ given by $g(x) = (1+x)^q - (1+qx)$. Then $f(0)=0$. The first-order derivative of $g$ is $g'(x) = q(1+x)^{q-1}-q$, so $g'(0)=0$. The second-order derivative of $g$ is $g''(x) = q(q-1)(1+x)^{q-2} \ge 0$ for all $x\in(-1,\infty)$ because $q\not\in(0,1)$. This means $g'$ is nonpositive on $(-1,0)$ and nonnegative on $(0,\infty)$, so $g$ attains its minimum at $0$. \Halmos\endproof}

\begin{lemma}\label{lemma:taylor-power-upper-bound}
For $x\ge-1/2$ and $q<0$, \[(1+x)^q \le 1+qx + 2^{-q+2}x^2.\]
\end{lemma}
{\proof{Proof of Lemma~\ref{lemma:taylor-power-upper-bound}} Denote the smooth function $f \in \mathcal{C}^\infty([-1,\infty))$ given by $f(x) = (1+x)^{q} - 1 - qx - 2^{-q+2} x^2$, then $f(-1/2) = -1+q/2 < 0$ and $f(0)=0$. The derivative of $f(x)$ is \[
f'(x) = q (1+x)^{q-1} - q - 2^{-q+3} x, 
\]
so $f'(0)=0$. Next, we argue that $0$ is the only local maximum of $f$.
The second-order derivative of $f$ is \[
f''(x) = q(q-1) (1+x)^{q-2} - 2^{-q+3}.\]
Note it is monotonically decreasing in $(-1, \infty)$ and $f''(0) = q(q-1)-2^{-q+3} < 0$. Now $f'(x) = f'(0) + \int_0^x f''(z)\ud z$ for any $x\in[-1, \infty)$, so for any $x\ne0$ where $f'(x)=0$, $f''(x) > 0$. This means it cannot be a local maximum. Hence, $0$ is the only local maximum of $f$.

Finally, notice $f'(0)=0$ and $f''(x)<0$ for any $x\ge0$ implies that $f(x)<f(0)=0$ for any $x\ge0$. Since we know $f(0)=0$, we conclude that $0$ is a global minimum of $f$ and that $f(x)\ge0$ for all $x\ge-1$. \Halmos \endproof}

\begin{lemma}
\label{lemma:series}
Given $\eta>1$ and $\zeta\in(0,1)$, for any $C>0$, \[ \sum_{k=1}^\infty \mathrm{e}^{\eta^k} \exp\left(-C\mathrm{e}^{\eta^k\zeta}\right) \le C^{-1/\zeta}h(\eta,\zeta), \]
where
\[ h(\eta,\zeta) := \q{1}{\zeta\log\eta}\Gamma\left(\q{1}{\zeta}\right) + 2\zeta^{-1/\zeta}.\]
\end{lemma}
\proof{Proof of Lemma~\ref{lemma:series}}
For ease of notation, we define a function $a:\RR\to\RR$ given by $a(w) = \eta^w -C\mathrm{e}^{\eta^w\zeta}$ so that the LHS of the inequality we are to show is $\sum_{k=1}^\infty \mathrm{e}^{a(k)}$. Clearly $a(\cdot)$ is continuously differentiable in $(0,\infty)$, and its first-order derivative is
\[a'(w) =  \left( 1- C\zeta \mathrm{e}^{\eta^w\zeta} \right) \eta^w \log\eta.\]
If $C\zeta\ge1$, then $a'(w)<0$ for all $w>0$, so that $a(\cdot)$ is monotonically decreasing, and $\sum_{k=1}^\infty \mathrm{e}^{a(k)} \le \int_0^\infty \mathrm{e}^{a(w)}\ud w$. Otherwise, $a'(w)>0$ and $a(w)$ is increasing for $w\in(0,w^*)$, and $a'(w)<0$ and $a(w)$ is decreasing for $w\in(w^*,\infty)$, where 
\[w^* := \q1{\log\eta} \log \left(\q{-\log(C\zeta)}{\eta}\right)>0.\]
This means $a(\cdot)$ is maximized at $w^*$, where it takes value 
\[a(w^*) = \q{-\log(C\zeta)-1}{\zeta}.\] 
Denoting $k^*:=\lfloor w^* \rfloor\ge0$, we have
\begin{equation*}
\begin{split}
    \sum_{k=1}^\infty \mathrm{e}^{a(k)} &\le \sum_{k=1}^{k^*-1} \mathrm{e}^{a(k)} + \e^{a(k^*)} + \e^{a(k^*+1)} + \sum_{k=k^*+2}^\infty \mathrm{e}^{a(k)} \\
    &\le \int_0^{k^*} \e^{a(w)}\ud w + 2\e^{a(w^*)} + \int_{k^*+1}^\infty \e^{a(w)}\ud w\\
    &\le \int_0^\infty \e^{a(w)}\ud w + 2(C\zeta)^{-1/\zeta},
\end{split}
\end{equation*}
where in the first inequality we assume $\sum_{k=p}^{q}\cdot = 0$ for any $p>q$, the second follows because $a(w)$ increases in $(0,k^*)$ and decreases in $(k^*+1, \infty)$, and the last by definition of $w^*$.

Now it remains to bound the integral $\int_0^\infty \e^{a(w)}\ud w$. In fact, we conclude using change of variables that
\begin{equation*}
\begin{split}
    \int_0^\infty \e^{a(w)}\ud w 
    &= \int_0^\infty \exp\left(\eta^w -C\mathrm{e}^{\eta^w\zeta}\right)\ud w \\
    &= \int_1^\infty \exp\left(u - C\e^{\zeta u}\right)\q{1}{u\log\eta}\ud u\\
    &\le \q1{\log\eta} \int_1^\infty \exp\left(u - C\e^{\zeta u}\right) \ud u \\
    &\le \q1{\log\eta} \int_0^\infty \exp\left(u - C\e^{\zeta u}\right) \ud u \\
    &= \q{C^{-1/\zeta}}{\log\eta} \int_0^\infty \exp\left(-v^\zeta\right)\ud v \\
    &= C^{-1/\zeta} \q{\Gamma(1/\zeta)}{\zeta\log\eta},
\end{split}
\end{equation*}
where the first equality is by definition of $a(\cdot)$, the second equality by change of variables $u = \eta^w$, the first inequality because $u\ge1$, the second inequality by non-negativity, the third equality by change of variables $v=C^{-1/\zeta}\e^u$, and the last equality by definition of the Gamma function $\Gamma(\cdot)$.\Halmos\endproof

\clearpage

\section{Additional policies}
\label{sec:add-alg}

\begin{algorithm}[htbp]
\caption{Backward Infrequent Re-solving ($\mathsf{BIR}_\eta$)}
\label{alg:bir}
\KwInput{time horizon length $T\in\N$, initial utilities $B_0\in\RR^n$, hyper-parameter $\eta\in(1, \nicefrac43)$.}
\KwInitialize{set $K \gets \lceil \log\log T / \log \eta \rceil$\;
set $t_0^*\gets0$, $t_k^* \gets T - \lfloor \exp(\eta^{K-k}) \rfloor$ for $k=1,\ldots,K$, and $t_{K+1}^* \gets T$\;}
\For{$k=0,1,\ldots,K$}{
    solve the fluid problem with updated utilities: \[\boldsymbol\xi^\F \in \arg\max \left\{ w \left( B_{t_k^*} + (T-t_k^*) \sum_{\ell\in[L]} p_\ell \beta_\ell  * \xi_\ell \right) : \boldsymbol\xi\in\Delta_n^L \right\}\]
    \For(\tcp*[h]{implement the updated fluid policy}){$t=t_k^*+1,\ldots,t_{k+1}^*$}{
        observe incoming supply with utility vector $b_t$\;
        \If{$b_t = \beta_\ell$}{
            act $x_t \gets (\xi^\F)_\ell$\;
        }
        update $B_t^i \gets B_{t-1}^i + b_t^i x_t^i$ for every agent $i \in [n]$\;
    }
}
\KwOutput{actions $(x_t:t=1,2,\ldots,T)$ and objective $\min_{i\in[n]} B_T^{i}$.}
\end{algorithm}

\begin{algorithm}[htbp]
\caption{Frequent Re-solving ($\mathsf{FR}$)}
\label{alg:fr}
\KwInput{time horizon length $T\in\N$, initial utilities $B_0\in\RR^n$.}
\KwInitialize{}
\For{$t=1,\ldots,T$}{
    solve the fluid problem with updated utilities: \[\boldsymbol\xi^\F \in \arg\max \left\{ w \left( B_{t-1} + (T-t+1) \sum_{\ell\in[L]} p_\ell \beta_\ell * \xi_\ell \right) : \boldsymbol\xi\in\Delta_n^L \right\}\;\]
    observe incoming supply with utility vector $b_t$\;
    \If{$b_t = \beta_\ell$}{
            act $x_t \gets (\xi^\F)_\ell$\;
        }
    update $B_t^i \gets B_{t-1}^i + b_t^i x_t^i$ for every agent $i \in [n]$\;
}
\KwOutput{actions $(x_t:t=1,2,\ldots,T)$ and objective $\min_{i\in[n]} B_T^{i}$.}
\end{algorithm}

\end{document}